\documentclass[12pt, reqno]{amsart}
\usepackage[utf8]{inputenc}
\usepackage[english]{babel}

\usepackage{ifthen}
\usepackage{microtype}
\usepackage[dvipsnames]{xcolor}
\usepackage{comment}
\usepackage{blkarray}
\usepackage{fullpage}
\usepackage{mathtools}
\usepackage{xfrac}
\usepackage{centernot}
\usepackage{enumitem}
\usepackage{mathdots}
\usepackage{multirow}
\usepackage{amssymb}
\usepackage{soul}
\usepackage{nicematrix}
\usepackage[noadjust]{cite}
\usepackage{booktabs}

\newlist{paraenum}{enumerate}{1}
\setlist[paraenum]{wide, label=(\arabic*)}

\makeatletter
\patchcmd{\@setaddresses}{\indent}{\noindent}{}{}
\patchcmd{\@setaddresses}{\indent}{\noindent}{}{}
\patchcmd{\@setaddresses}{\indent}{\noindent}{}{}
\patchcmd{\@setaddresses}{\indent}{\noindent}{}{}
\makeatother

\usepackage{caption}
\usepackage{subcaption}
\usepackage{tikz}    
\usetikzlibrary{cd}
\usetikzlibrary{backgrounds}
\usetikzlibrary{shapes}
\usetikzlibrary{arrows}
\usetikzlibrary{positioning}
\usetikzlibrary{calc}
\usetikzlibrary{quotes}
\usetikzlibrary{fit}
\usetikzlibrary{decorations.pathreplacing}
\definecolor{cimblue}{RGB}{38,70,150}
\definecolor{cimpurple}{RGB}{122,58,170}
\colorlet{cimedge}{cimblue!85!black}
\colorlet{cimnode}{cimblue!80!black}
\colorlet{cimnodefill}{cimblue!7}
\colorlet{cimcoll}{cimpurple!25}
\colorlet{cimhl}{cimpurple!90!black}
\tikzset{
  coll/.style={fill=cimcoll,draw=cimpurple!85!black},
  vtx/.style={circle,draw=cimnode,fill=cimnodefill,inner sep=0pt,minimum size=12pt,font=\scriptsize},
  tedge/.style={-{Stealth[length=3.5pt]},thick,cimedge},
  xedge/.style={-{Stealth[length=3.5pt]},dashed,thick,cimhl},
  hvtx/.style={draw=cimnode,fill=cimnodefill,rounded corners=2pt,inner sep=2pt,minimum width=14pt},
  harc/.style={-{Stealth[length=3.5pt]},cimedge},
}

\usepackage[pdftex]{hyperref}
\hypersetup{
    colorlinks=true,
    linkcolor={magenta},
    citecolor={blue},
    urlcolor={blue}
}
\usepackage[nameinlink,capitalise]{cleveref}
\numberwithin{equation}{section}
\theoremstyle{plain}
\newtheorem{theorem}{Theorem}[section]

\newtheorem{corollary}[theorem]{Corollary}
\newtheorem{lemma}[theorem]{Lemma}
\newtheorem{proposition}[theorem]{Proposition}

\theoremstyle{definition}
\newtheorem{remark}[theorem]{Remark}
\newtheorem{examplex}[theorem]{Example}
\newenvironment{example}
  {\pushQED{\qed}\examplex}
  {\popQED\endexamplex}

\newtheorem{definition}[theorem]{Definition}

\newcommand{\rr}{\mathbb{R}}

\newcommand{\cb}{\mathcal{B}}

\newcommand{\cg}{\mathcal{G}}

\newcommand{\cv}{\mathcal{V}}

\newcommand{\ch}{\mathcal{H}}
\newcommand{\Nb}{\mathcal N}

\newcommand{\ba}{\mathbf{a}}
\newcommand{\bc}{\mathbf{c}}
\newcommand{\bp}{\mathbf{p}}
\newcommand{\bx}{\mathbf{x}}
\newcommand{\be}{\mathbf{e}}

\newcommand{\skel}{\mathrm{skel}}

\DeclareMathOperator{\rank}{rank}

\newcommand{\cim}{\mathrm{CIM}}
\newcommand{\bdcim}[2]{\mathrm{CIM}_{#1, #2}}
\newcommand{\dags}[2]{\mathfrak{G}_{#1,#2}}

\newcommand{\pa}{\mathrm{pa}}
\newcommand{\conv}{\mathrm{conv}}
\newcommand{\coll}{\mathrm{coll}}
\renewcommand{\rightarrow}{\!\to\!}

\definecolor{benpurple}{RGB}{180, 0, 240}

\title{Polyhedral Geometry of Degree-Bounded Characteristic Imset Polytopes}

\author{Jane Ivy Coons}
\address{Worcester Polytechnic Institute}
\email{jcoons@wpi.edu}

\author{Benjamin Hollering}
\address{Max Planck Institute for Mathematics in the Sciences, Inselstra{\ss}e 22, 04103 Leipzig, Germany}
\email{benjamin.hollering@mis.mpg.de}

\author{Niharika Chakrabarty Paul}
\address{Max Planck Institute for Mathematics in the Sciences, Inselstra{\ss}e 22, 04103 Leipzig, Germany}
\email{niharika.paul@mis.mpg.de}

\begin{document}

\begin{abstract}
    Characteristic imset polytopes are convex hulls of a collection of characteristic imset vectors, which are 0/1 vectors that encode the Markov equivalence class of a directed acyclic graph. These polytopes are the feasible region of linear programs which correspond to causal discovery problems. We introduce a new family of these polytopes whose vertices correspond to directed acyclic graphs with in-degree bounded by $k$. This yields a polytope whose dimension is only polynomial in the number of nodes of the underlying graphs instead of exponential. We show that when $k = 1$, this polytope is a matroid independence polytope. For $k = 2$ we provide several large families of facet-defining inequalities, which give a complete $\mathcal{H}$-representation of the polytope on four nodes.
\end{abstract}

\maketitle

\section{Introduction}
\noindent Directed acyclic graphical (DAG) models use a DAG $\mathcal{G} = ([n], E)$ to model the causal relationships between events which correspond to the components of a random vector $X = (X_1, \ldots, X_n)$. The associated statistical model $M_\mathcal{G}$ is the set of all distributions which are \emph{Markov} to the DAG, meaning that the joint density admits a factorization
\[
f(x_1, \ldots, x_n) = \prod_{i = 1}^n f(x_i | x_{\pa(i)}),
\]
where $\pa(i)$ is the set of \emph{parents} of the node $i$ in $\mathcal{G}$. DAG models are central objects in causality and have found many applications throughout the sciences including biology, medicine, and economics \cite{friedman2000using,sachs2005causal,greenland1999causal,imbens2020potential}. In many of these applications, one has access to data which are believed to be generated by a distribution which belongs to some true but unknown DAG model $M_\mathcal{G}$, and would like to reconstruct the graph $\cg$ from this data. This problem is often known as \emph{causal discovery} and is a central problem in causality. 

An inherent difficulty when studying these models is that distinct DAGs $\cg$ and $\ch$ may produce the same set of distributions, that is, $M_\cg = M_\ch$. Such DAGs are called \emph{Markov equivalent}, and there exists a well-known combinatorial characterization of the equivalence classes \cite{verma1990equivalence,andersson1997characterization}. This means that only the equivalence class of the true DAG can be reconstructed, given access to only observational data. Despite this, there exist many causal discovery algorithms which are often divided into the following categories. \emph{Constraint-based} causal discovery algorithms, such as the well-known \emph{PC algorithm}, find a graph such that all distributions in the model satisfy the same constraints, for example conditional independencies, as those exhibited by the given data \cite{spirtes1991algorithm,spirtes2000causation}. On the other hand, \emph{score-based} algorithms such as \emph{Greedy Equivalence Search} optimize a score function such as the Bayesian Information Criterion (BIC) as they iteratively mutate the DAG via edge additions, reversals, and deletions \cite{chickering2002optimal}.

The introduction of \emph{characteristic imsets} \cite{studeny2010characteristic,hemmecke2012characteristic} led to many new algorithms and perspectives on causal discovery. The characteristic imset vector of a DAG $\cg = ([n], E)$ is a $0/1$ vector $\bc_\cg \in \rr^{2^n - n -1}$ which encodes the Markov equivalence class of a DAG. That is, $\bc_\cg = \bc_\ch$ if and only if $\cg$ and $\ch$ are Markov equivalent. Moreover, Studen\'y et al. \cite{studeny2010characteristic} showed that the problem of causal discovery can be seen as a linear program over the polytope $\cim_{n} = \conv(\bc_\cg : \cg \text{ DAG on } [n])$. In particular, the BIC can be seen as a vector $\boldsymbol\alpha \in \rr^{2^n - n -1}$ and thus the DAG which best fits given data can be found by solving the linear program
\begin{align*} 
    \max_{\bx \in \cim_n} \boldsymbol\alpha \cdot \bx.
\end{align*}
This new geometric perspective on causal discovery has led to a much better understanding of score-based causal discovery algorithms, and yielded improved greedy algorithms such as \emph{skeletal greedy CIM} \cite{linusson2023greedy, linusson2026edges}. Unfortunately, performing causal discovery by solving this linear program explicitly still remains quite challenging for a variety of reasons. First, the polytope $\cim_n$ lives in a space which grows exponentially in $n$. Thus, even if a polynomial time solver can be employed to solve the above linear program, the resulting algorithm is still exponential in the number $n$ of nodes of the underlying DAG. Moreover, the most efficient algorithms for solving linear programs, such as interior point methods, require that the feasible region be presented as an intersection of half-spaces, but this description is generally not known for $\cim_n$. Some recent progress has been made on this front by the second author together with Johnson and Solus \cite{johnson2026hyperplane} as well as by Xi and Yoshida \cite{xi2015characteristic}; however, their work focuses on a smaller subpolytope of $\cim_n$ where only graphs with a fixed skeleton are considered in the convex hull or a total order on the vertices is assumed to be known.

In this paper, we suggest an alternative to $\cim_n$ which is naturally motivated by causal discovery and helps combat the shortfalls of the above linear programming approach. A common assumption by many classic causal discovery algorithms is that the maximum in-degree of any vertex in the true graph is bounded by some fixed integer $k$. For example, the PC algorithm is only polynomial time under this assumption, since it would otherwise need to check every possible pairwise conditional independence statement of which there are $\binom{n}{2}2^{n-2}$. This naturally leads us to consider the \emph{degree-bounded} CIM polytope
\[
\cim_{n, k} := \conv(\bc_\cg ~:~ \cg \text{ a DAG such that }\max_{i \in [n]} |\pa_\cg(i)| \leq k). 
\]
Unlike $\cim_n$, the degree-bounded CIM polytope has dimension $O\left(\binom{n}{k+1}\right) = O(n^{k+1})$ which grows only polynomially in $n$ for fixed $k$. Our main objective in this paper is thus focused on determining an \emph{H-representation} of $\cim_{n, k}$, particularly in the case where $k = 1, 2$ since this would aid in the efficient estimation of graphs with small in-degree via the linear programming approach described above. Note that $k = 1, 2$ includes commonly used classes of graphs such as binary trees or phylogenetic networks. 

We obtain a complete description of these inequalities for $\cim_{n, 1}$ by showing that it is in fact a \emph{matroid independence polytope}. We then describe several families of facet-defining inequalities for $\cim_{n,2}$, which give a complete $\ch$-representation of $\cim_{4,2}$ together with some obvious inequalities. More specifically, we show that every connected undirected graph yields an inequality, which we call the \emph{path inequality}, while every undirected tree without vertices of degree $2$ yields an additional inequality, which we call the \emph{outer inequality}. This suggests that despite the fact that $\cim_{n,2}$ is much smaller than $\cim_n$, the number and complexity of these inequalities still grow quickly with $n$.

The remainder of this paper is organized as follows. In Section \ref{sec:prelim}, we provide some background information on polytopes and directed acyclic graphs and their models. In Section \ref{sec:deg-bounded-cim}, we introduce the degree-bounded characteristic imset polytope, compute its dimension, and study its relationship to $\cim_n$. We then quickly show that $\cim_{n,1}$ is a matroid independence polytope and thus obtain a complete $\mathcal{H}$-description of it. Sections \ref{sec:path-ineqs}, \ref{sec:star-tree-ineqs}, and \ref{sec:clique-ineqs} introduce families of inequalities for $\cim_{n, 2}$ which are naturally associated to subgraph structures such as connected graphs, trees, and cliques. In each case, we show that all of the stated inequalities are valid and facet-defining. We end with a discussion of the completeness of our inequalities and directions for future work.

\section{Preliminaries}
\label{sec:prelim}
\subsection{Polytope Basics}
\noindent We will first provide some background on convex polytopes that we will be using throughout the paper. We refer the reader to \cite{ziegler2012lectures} for a more detailed exposition. 
\begin{definition}
    An \emph{$\mathcal{H}$-polyhedron} is ${P}=\{\bx\in\mathbb{R}^d\vert A\bx \leq \mathbf{b}\}$, where $A\in \mathbb{R}^{n\times d}$ and $\mathbf{b}\in\mathbb{R}^n$. In other words, an $\mathcal{H}$-polyhedron is an intersection of a finite number of half-spaces. Furthermore, an \emph{$\mathcal{H}$-polytope} is a bounded $\mathcal{H}$-polyhedron. 

    A \emph{$\mathcal{V}$-polytope} is defined as the convex hull of a finite set of vectors. Here, the convex hull of a set of points $\bp_1,\dots,\bp_\ell$ is defined as follows. 
    $$\mathrm{conv}(\bp_1,\dots,\bp_\ell)=\left\{ \sum_{i\in[\ell]} \lambda_i\bp_i \Bigg\vert\ \lambda_i\geq 0,\ \sum_{i\in[\ell]} \lambda_i = 1\right\}.$$
\end{definition}

\noindent The following theorem describes the connection between the two categories of polytopes, as defined above. 

\begin{theorem}[{\cite[Theorem 1.1]{ziegler2012lectures}}]
    A subset ${P}\subset\mathbb{R}^d$ is a $\mathcal{V}$-polytope if and only if it is an $\mathcal{H}$-polytope.
\end{theorem}

\noindent The proof of this theorem is algorithmic, via a process called \emph{Fourier-Motzkin elimination}, which computes an $\mathcal H$-representation of a polytope from a $\mathcal V$-representation, and vice versa \cite[Lecture~1]{ziegler2012lectures}. In principle, one can, therefore, always pass from one description to the other, although the number of inequalities produced can grow very quickly with the dimension. In this paper, we are interested in describing the $\ch$-representation of a special class of polytopes. This is equivalent to listing its \emph{facets}. We define this term, and related concepts, below. 

\noindent Throughout, for $\ba,\bx\in\rr^d$ we write $\ba\cdot\bx=\sum_{i}a_ix_i$ for the standard inner product. When a linear functional on $\rr^d$ is given a name, such as $\ell$, we write $\ell(\bx)$ for its value at $\bx$; thus $\ell(\bx)=\ba\cdot\bx$, where $\ba$ is the vector of coefficients of $\ell$.

\begin{definition}\label{def: hat}
    Let $P\subset\mathbb{R}^d$ be a convex polytope. A linear inequality $\ba\cdot\bx\leq a_0$, with $\ba\in\rr^d$ and $a_0\in\rr$, is \emph{valid} on the polytope if it holds for all points $\bx$ in the polytope. A \emph{face} of $P$ is a set $F$, given by:
    $$F = P\cap \{\bx\in\mathbb{R}^d\mid\ba\cdot\bx= a_0\},$$
    where $\ba\cdot\bx\leq a_0$ is a valid inequality on $P$. We define the dimension of $F$ to be the dimension of its affine hull. The dimension $0$ faces of a polytope are the \emph{vertices}, the dimension $1$ faces the \emph{edges}, and the faces of dimension $\text{dim}(P)-1$ are the \emph{facets} of the polytope. 
\end{definition}

\noindent In order to verify whether or not the valid inequality $\ba\cdot\bx\leq a_0$ does indeed define a facet on a polytope $P$, we can undertake the following steps. 
\begin{enumerate}
    \item We write down the matrix $M$ such that each column of $M$ is a vertex of the face $F = P\cap \{\bx\in\mathbb{R}^d\mid\ba\cdot\bx= a_0\}$. We will sometimes refer to this matrix as the \emph{defining matrix} of the face. 
    \item Then, we ensure that the vector of ones is contained in the rowspan of $M$. If it is not, we augment the matrix with the row of ones and call this new matrix $\widehat{M}$. If the vector of ones already belongs to the rowspan of $M$, we set $\widehat{M}$ as $M$. 
    \item Then, the dimension of the affine hull of $F$ is one less than $\text{rank}(\widehat{M})$. 
    \item Therefore, $F$ is a facet of $P$ if and only if $\text{rank}(\widehat{M}) = \text{dim}(P)$.
\end{enumerate}
We illustrate this algorithm through the following simple example. 
\begin{example}
    Consider the cube in three dimensions. This is a convex polytope with $\cv$-representation given by the columns of the following matrix:
    $$\begin{pmatrix}
        0 & 1 & 0 & 0 & 1 & 1 & 0 & 1\\
        0 & 0 & 1 & 0 & 1 & 0 & 1 & 1\\
        0 & 0 & 0 & 1 & 0 & 1 & 1 & 1
    \end{pmatrix}$$
    This polytope has 8 vertices, 12 edges, and 6 facets. One such facet is specified by the inequality $(0,0,-1)\cdot\bx\leq 0$. This is the bottom-most face, $F$. Then, we can write down the defining matrix of this face as $M = \begin{pmatrix}
        0 & 1 & 0 & 1\\0 & 0 & 1 & 1\\0 & 0 & 0 & 0
    \end{pmatrix}$. Clearly, $(1~1~1~1)$ is not in the rowspan of $M$. Thus, we define the matrix $\widehat{M}$ as:
    $$\widehat{M} = \begin{pmatrix} 1& 1&1&1\\
        0 & 1 & 0 & 1\\0 & 0 & 1 & 1\\0 & 0 & 0 & 0 \end{pmatrix}.$$ 
    We can then calculate the rank of $\widehat{M}$ as $3 = \text{dim}(P)$. This calculation shows that $F$ is indeed a facet. 
\end{example}

\subsection{Directed Acyclic Graphical Models}
In this subsection, we set notation concerning structures in directed graphs and provide a brief refresher on graphical models on DAGs. For a more detailed introduction to this topic, we refer the reader to \cite{lauritzen1996graphical,pearl2009causality,sullivant2018algebraic}.  
Graphical models use graphs to represent conditional independence and causal relationships between random variables. A density $f(x_1, \ldots, x_n)$ belongs to the graphical model on a DAG $\cg = ([n], E)$ if it admits a factorization of the form
\[
f(x_1, \ldots, x_n) = \prod_{i = 1}^n f(x_i | x_{\pa(i)}),
\]
where $\pa(i)$ is the set of parents of the node $i$ in $\cg$. 

As discussed previously, different DAGs do not always yield different sets of distributions. Two DAGs are said to be \emph{Markov equivalent} if the sets of densities that factorize according to the two graphs are the same. The following celebrated theorem of Verma and Pearl gives a graphical characterization of Markov equivalence. In this theorem, the \emph{skeleton} of $\cg$ refers to the underlying undirected graph of $\cg$ which we denote by $\skel(\cg)$. A \emph{collider} in a DAG is a path of the form $i\rightarrow j \leftarrow k$. It is \emph{unshielded} if $i$ and $k$ are not adjacent, that is, if the path is induced, and the set of all unshielded colliders is denoted by $\coll(\cg)$. We note that unshielded colliders are also often called \emph{v-structures} in the literature.

\begin{theorem}[{\cite{verma1990equivalence}}]
\label{thm:meq}
    Two DAGs $\cg$ and $\ch$ are Markov equivalent if and only if they have the same skeleton and unshielded colliders. 
\end{theorem}

Unshielded colliders can also be described through \emph{moralization}. The \emph{moral graph} $\cg^m$ of a DAG $\cg$ is the undirected graph obtained from $\skel(\cg)$ by adding an edge between every two vertices that have a common child. An edge is added between $i$ and $k$ exactly when $i\to j\leftarrow k$ for some $j$ and $i,k$ are not adjacent in $\cg$, that is, when $i\to j\leftarrow k$ is an unshielded collider \cite{lauritzen1996graphical}. For $A\subseteq[n]$ let $\cg[A]$ be the subgraph of $\cg$ induced on $A$. This gives the following reformulation of Theorem~\ref{thm:meq}.
 
\begin{corollary}\label{cor:moral}
    Two DAGs $\cg$ and $\ch$ on $[n]$ are Markov equivalent if and only if $(\cg[A])^m=(\ch[A])^m$ for every $A\subseteq[n]$. It suffices to check this for all $A$ with $|A|\le 3$.
\end{corollary}
\begin{proof}
    The moral graph of a DAG is determined by its skeleton and its unshielded colliders, and the skeleton and unshielded colliders of $\cg[A]$ are those of $\cg$ that lie inside $A$. So Markov equivalent DAGs satisfy $(\cg[A])^m=(\ch[A])^m$ for all $A$. Conversely, if this holds for all $A$ with $|A|\le3$, then the sets $A$ of size $2$ show that $\skel(\cg)=\skel(\ch)$. For a set $A=\{i,j,k\}$ with $i\sim j\sim k$ and $i\not\sim k$, the moral graph of $\cg[A]$ contains the edge $ik$ if and only if $i\to j\leftarrow k$ is an unshielded collider of $\cg$. So $\coll(\cg)=\coll(\ch)$, and Theorem~\ref{thm:meq} applies.
\end{proof}
 
\noindent Note that it is not enough to compare the moral graphs of $\cg$ and $\ch$ themselves when testing for Markov equivalence. For instance, the DAGs $1\to 2\leftarrow 3\leftarrow 4\to 1$ and $2\to1\to4\leftarrow3\leftarrow2$ have the same skeleton-the undirected $4$-cycle $1-2-3-4-1$-and the same moral graph, but their unshielded colliders are at $2$ and at $4$ respectively. We now look at the following example of two DAGs in the same Markov equivalence class.

\begin{example}\label{ex:mec}
    Consider the following two DAGs on $[4]$.
    
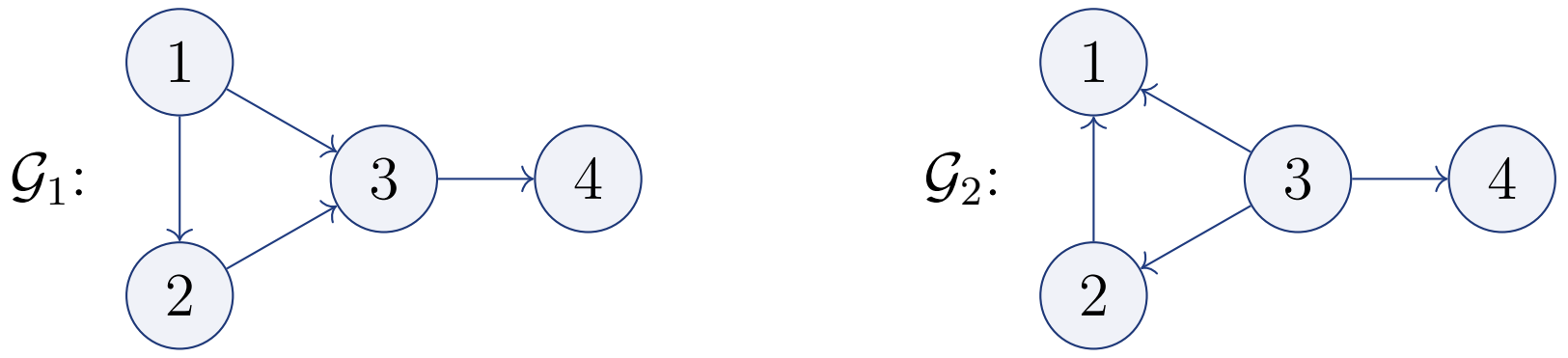
\begin{center}
    \begin{tikzpicture}[every node/.style={minimum size=0.5cm}]
    \node at (-0.9,0) {$\cg_1$:};
    \node [shape=circle, draw=cimnode, fill=cimnodefill] (A) at (0,0.8) {$1$};
    \node [shape=circle, draw=cimnode, fill=cimnodefill] (B) at (0,-0.8) {$2$};
    \node [shape=circle, draw=cimnode, fill=cimnodefill] (C) at (1.4,0) {$3$};
    \node [shape=circle, draw=cimnode, fill=cimnodefill] (D) at (2.8,0) {$4$};
    \path [->,cimedge] (A) edge (B);
    \path [->,cimedge] (A) edge (C);
    \path [->,cimedge] (B) edge (C);
    \path [->,cimedge] (C) edge (D);
    \end{tikzpicture}
    \qquad\qquad
    \begin{tikzpicture}[every node/.style={minimum size=0.5cm}]
    \node at (-0.9,0) {$\cg_2$:};
    \node [shape=circle, draw=cimnode, fill=cimnodefill] (A) at (0,0.8) {$1$};
    \node [shape=circle, draw=cimnode, fill=cimnodefill] (B) at (0,-0.8) {$2$};
    \node [shape=circle, draw=cimnode, fill=cimnodefill] (C) at (1.4,0) {$3$};
    \node [shape=circle, draw=cimnode, fill=cimnodefill] (D) at (2.8,0) {$4$};
    \path [->,cimedge] (B) edge (A);
    \path [->,cimedge] (C) edge (A);
    \path [->,cimedge] (C) edge (B);
    \path [->,cimedge] (C) edge (D);
    \end{tikzpicture}
    \end{center}
    Both DAGs have the skeleton with edges $12,13,23,34$. In each of them exactly one triple of vertices has a collider, namely $\{1,2,3\}$, but the collider is at a different vertex: in $\cg_1$ we have $1\to3\leftarrow2$, while in $\cg_2$ we have $2\to1\leftarrow3$. By Theorem~\ref{thm:meq}, $\cg_1$ and $\cg_2$ are Markov equivalent.
\end{example}

\section{Degree-Bounded CIM Polytopes}
\label{sec:deg-bounded-cim}
\noindent In this section, we introduce the degree-bounded characteristic imset (CIM) polytope which is the core object of study in this paper. We first recall the classical CIM polytope of \cite{studeny2010characteristic}.

\begin{definition}
    Let $[n]:= \{1,\dots,n\}$ and $\cg$ be a directed acyclic graph (DAG), with vertex set $[n]$. Then, we define the characteristic imset $\bc_{\cg}$ as the function:
    \begin{align*}
        \bc_{\cg}:  \{S\subseteq [n] : |S|\geq 2\} &\longrightarrow ~ \{0,1\}\\
                 S &\longmapsto  ~ 
\begin{cases}
1 \text{ if there exists some } i \in S \text{ such that } S\setminus \{i\} \subseteq \mathrm{pa}_{\cg} (i)\\
0 \text{ otherwise}.
\end{cases}
    \end{align*}   
\noindent The characteristic imset polytope (CIM polytope), denoted $\cim_n$, is given by
$$\cim_n := \mathrm{conv}(\bc_{\cg} : \cg \text{ is a DAG on } [n]).$$
\end{definition}

The characteristic imset is of interest because it is a complete invariant of Markov equivalence.
 
\begin{theorem}[{\cite{studeny2010characteristic,hemmecke2012characteristic}}]\label{thm:cim-mec}
    Two DAGs $\cg$ and $\ch$ on $[n]$ are Markov equivalent if and only if $\bc_\cg=\bc_\ch$.
\end{theorem}
 
In fact, only the coordinates of $\bc_\cg$ indexed by sets of size $2$ and $3$ are needed. Indeed, $\bc_\cg(\{i,j\})=1$ if and only if $i$ and $j$ are adjacent, so these coordinates record $\skel(\cg)$. For a $3$-set $\{i,j,k\}$ with $i\sim j\sim k$ and $i\not\sim k$, we have $\bc_\cg(\{i,j,k\})=1$ if and only if $i\to j\leftarrow k$, so these coordinates record $\coll(\cg)$. By Theorem~\ref{thm:meq}, the coordinates of size $2$ and $3$ therefore determine the Markov equivalence class of $\cg$, and by Theorem~\ref{thm:cim-mec} they determine all of $\bc_\cg$.

 We now illustrate the characteristic imset on the DAGs of Example~\ref{ex:mec}.

\begin{example}
For a DAG $\cg$, the value $\bc_{\cg}(\{i,j\})$ is $1$ if and only if $i$ and $j$ are adjacent, that is, there is an edge between $i$ and $j$. The value $\bc_{\cg}(\{i,j,k\})$ is $1$ if and only if one of the three vertices, say $i$, has the other two as parents, i.e.\ $j\to i\leftarrow k$. Similarly, $\bc_\cg([4])=1$ if and only if some vertex has the other three as parents. 
 
For $\cg_1$ of Example~\ref{ex:mec}, the parent sets are $\pa(1)=\emptyset$, $\pa(2)=\{1\}$, $\pa(3)=\{1,2\}$ and $\pa(4)=\{3\}$. For $\cg_2$ they are $\pa(1)=\{2,3\}$, $\pa(2)=\{3\}$, $\pa(3)=\emptyset$ and $\pa(4)=\{3\}$. In both DAGs, no vertex has three parents, so $c([4])=0$, and the only $3$-set containing a vertex together with two of its parents is $\{1,2,3\}$: in $\cg_1$ this vertex is $3$, in $\cg_2$ it is $1$. Hence, the two characteristic imsets, written as row vectors, are as follows. 
\[\begin{blockarray}{cccccccccccc}
 & 12 & 13 & 14 & 23 & 24 & 34 & 123 & 124 & 134 & 234 & 1234\\
\begin{block}{c(ccccccccccc)}
\bc_{\cg_1} & 1 & 1 & 0 & 1 & 0 & 1 & 1 & 0 & 0 & 0 & 0\\
\bc_{\cg_2} & 1 & 1 & 0 & 1 & 0 & 1 & 1 & 0 & 0 & 0 & 0\\
\end{block}
\end{blockarray}.\]
Here, each column is indexed by the subset of $[4]$ given above it. So, $\bc_{\cg_1}=\bc_{\cg_2}$, as predicted by Theorem~\ref{thm:cim-mec}. The coordinate $123$ only records that $\{1,2,3\}$ has a collider, not at which vertex, so it does not distinguish the two DAGs.
\end{example}

Motivated by causal discovery, there has been a growing interest in the geometric structure of $\cim_n$. This includes a description of some of its edges and their relationship to greedy equivalence search \cite{linusson2023greedy,linusson2026edges}, as well as several papers which focus on computing the facet description of certain subpolytopes or faces of $\cim_n$ \cite{studeny2013polyhedral,xi2015characteristic,cussens2017polyhedral,johnson2026hyperplane}. Most work so far has focused on the faces $\cim_G = \conv(\bc_\cg ~:~ \skel(\cg) = G)$, since they have considerably fewer vertices and their dimension is much smaller. An understanding of their facets then enables so-called \emph{hybrid} causal discovery algorithms where one first estimates a skeleton $G$ and then orients the edges via solving a linear program over the smaller polytope $\cim_G$. This work, though, is still limited to the case where the underlying skeleton is a tree.

In this paper we study a new subpolytope of $\cim_n$ which is also naturally motivated by causal discovery. Recall that many classic causal discovery algorithms require a bound on the degree, or some similar sparsity assumption, to guarantee that the runtime of the algorithm is polynomial. For instance, the PC algorithm is essentially $O(n^{k+1})$ where $n$ is the number of nodes in the graph and $k$ is the maximum in-degree of any node in the graph. This suggests the following definition.

\begin{definition}
    Let $\dags{n}{k}$ denote the set of DAGs on $[n]$ in which every vertex has at most $k$ parents, that is, $|\pa_\cg(i)|\leq k$ for all $i\in[n]$. The \emph{degree-bounded CIM polytope} on $[n]$ with degree bound $k$ is defined as
    $$\bdcim{n}{k} := \mathrm{conv}(\bc_{\cg} : \cg \in \dags{n}{k}).$$
\end{definition}

$\bdcim{n}{k}$ has many immediate advantages over the full polytope $\cim_n$. First, note that if $\bc_\cg$ is a vertex 
of $\cim_{n, k}$, then $\bc_\cg(S) = 0$ for any subset $S \subset[n]$ such that $|S| > k+1$. Thus, $\bdcim{n}{k}$ naturally lives in a space of dimension $O(\binom{n}{k+1}) = O(n^{k+1})$ and thus the dimension of the ambient space grows only polynomially in the number of nodes $n$ for fixed $k$. If we have an explicit $\ch$-representation of $\cim_{n, k}$, then we immediately get an algorithm for causal discovery by writing the BIC as a linear functional in $\rr^{2^n-n-1}$ and optimizing it over $\cim_{n, k}$, under the assumption that the true graph has bounded in-degree. Since the complexity of solving an LP via interior point methods is polynomial in the number of constraints and number of decision variables, this resulting causal discovery algorithm will be polynomial in the number of facets of $\bdcim{n}{k}$ and the number of nodes $n$. So as in previous work on the CIM polytope and its subpolytopes, a natural first question is to understand the facets of $\cim_{n, k}$, which will be the primary focus of the remainder of this paper. Our first result shows that $\bdcim{n}{k}$ is full-dimensional in its ambient space. 

\begin{proposition}
$\bdcim{n}{k}$ is a face of $\cim{n}$ and full-dimensional in its ambient space and thus $\dim \bdcim{n}{k} = \sum_{i = 2}^{k+1} \binom{n}{i}$. 
\end{proposition}
\begin{proof}
First, observe that 
\[
\bdcim{n}{k} = \cim_n \cap \{x \in \rr^{2^{n} - n - 1} ~:~ x_S = 0 \text{ for all } S \subset [n]  \text{ such that } |S| > k+1\}. 
\]
The latter set is a coordinate linear space and thus we see that $\bdcim{n}{k}$ is a face of $\cim_n$. To see that $\bdcim{n}{k}$ is full-dimensional, fix the standard total ordering on $[n]$ and consider the complete DAG $K_S$ on the node set $S \subset [n]$ with orientation $i \to j$ if $i < j$. For each $S$ such that $2 \leq |S| \leq k+1$, $K_S \in \mathfrak{G}_{n,k}$ and the affine span of the characteristic imsets of these DAGs $K_S$ has dimension $\dim \bdcim{n}{k} = \sum_{i = 2}^{k+1} \binom{n}{i}$. Since the zero vector lies in $\bdcim{n}{k}$, the second point follows from the fact that the matrix of $\bc_{K_S}$ with columns and rows in graded lexicographic order on the subsets $S$ is upper triangular with 1's on the diagonal. 
\end{proof}

We now turn to the study of $\bdcim{n}{k}$ for small $k$, and thus naturally begin with $\bdcim{n}{1}$. The vertices of $\bdcim{n}{1}$ are the characteristic imsets of DAGs with in-degree at most $1$, that is, of directed forests in which every connected component is an \emph{arborescence}. We now show that $\cim_{n, 1}$ is in fact a \emph{matroid independence polytope}. First, we recall the definition of a matroid.  

\begin{definition}\label{def:matroid}
A \emph{matroid} $M = (E,\mathfrak{I})$ is a pair where $E$ is a finite set and $\mathfrak{I} \subseteq 2^E$ satisfies
\begin{enumerate}
    \item $\emptyset \in \mathfrak{I}$.
    \item If $I' \subseteq I \in \mathfrak{I}$, 
    then $I' \in \mathfrak{I}$.
    \item If $I_1, I_2 \in \mathfrak{I}$ and $|I_2| > |I_1|$, then there exists $e \in I_2 \setminus I_1$ such that $I_1 \cup \{e\} \in \mathfrak{I}$.
\end{enumerate}
\end{definition}

Matroids are combinatorial objects which provide an abstraction of the notion of independence and have been studied before in the context of statistics and conditional independence \cite{matus1999matroid,studeny2020contribution}. Every matroid has several polytopes naturally associated to it, one of which is the matroid independence polytope
\[
P(M) = \text{conv}(\chi_I\vert I\in \mathfrak{I}).
\]
Here, for each $A\subseteq E$, $\chi_A$ is the characteristic or indicator vector of $A$ in $\mathbb{Z}^{\vert E \vert}$. The following well-known theorem gives a complete characterization of the facets of such polytopes in terms of the flats of the matroid, which are subsets of the ground set to which adding any further element increases the rank $r$.

\begin{theorem}[{\cite[Theorem 40.5]{schrijver2003combinatorial}}]\label{thm: indep matroid facets}
    Suppose that $M$ has no loops. Then, $\text{dim}(P(M)) = |E|$. Furthermore, $P(M)$ is completely described by the following list of facets.
    \begin{itemize}
        \item $x(e)\geq 0,$ for all $e\in E$.
        \item $\sum_{e\in F} x(e) \leq r(F)$, for each indecomposable flat $F$.
    \end{itemize}
\end{theorem}

One particularly nice family of matroids are \emph{graphic matroids}. Consider an undirected graph $G = (V, E)$. The graphic matroid of $G$ has ground set $E$, and a subset $I \subseteq E$ is independent if the graph $(V,I)$ is a forest, that is, it contains no cycle. The following result shows that $\cim_{n, 1}$ is the matroid independence polytope of the graphic matroid of $K_n$. 

\begin{proposition}\label{prop: matroid independence}
The degree-bounded characteristic imset polytope $\cim_{n,1}$ equals the matroid independence polytope $P(K_n)$ and thus is cut out by the inequalities
\begin{itemize}
    \item $x_{ij}\geq 0$, for all distinct $i,j\in [n]$,
    \item for each $k\in \{2,\ldots,n\}$ and each set $\{i_1,\ldots,i_k\} \subseteq [n]$ of size $k$,
    $$\sum_{1\leq a<b\leq k} x_{i_ai_b} \leq k-1.$$
\end{itemize}    
\end{proposition}

The facet description of $P(K_n)$ in Proposition \ref{prop: matroid independence} has appeared before in \cite{edmonds1971matroids}. However, we provide the following proof here for completeness.

\begin{proof}
Consider the graphic matroid $M$ on the complete graph $K_n$ so the ground set is $E = \{x_{ij}, i,j\in [n]\}$, which is our set of edge variables. Then the independent sets $\mathfrak{I}$ are simply all forests on $n$ nodes and thus $\bdcim{n}{1}$ is equal to $P(M)$.

It remains to compute the indecomposable flats associated to this polytope. For $F\subseteq E$, let $(V,F)$ denote the spanning subgraph of $K_n$ with edge set $F$, so that $V=[n]$. Then, $F$ is a flat of the matroid if and only if every edge $e\in E\setminus F$ has its endpoints in two different connected components of $(V,F)$. Equivalently, 
$$F=\{e\in E: \text{ both endpoints of } e \text{ lie in the same connected component of } (V,F)\}.$$
Furthermore, a flat $F$ is indecomposable if and only if $(V,F)$ has exactly one nontrivial connected component. That is, each indecomposable flat $F$ of $M$ is the edge set of the complete graph on some set $\{i_1,\ldots,i_k\} \subseteq [n]$ with $k\geq 2$. The rank of a flat $F$ is $|V|$ minus the number of connected components of $(V,F)$, which equals $k-1$ for this flat. Thus, we can rewrite the facets given in Theorem \ref{thm: indep matroid facets} as the list in the statement of this proposition. 
\end{proof}

 We now turn our attention to $\cim_{n, 2}$ which is our primary focus for the remainder of this paper. 
As we will show in the later sections, $\cim_{n, 2}$ already has an extremely rich and complex combinatorial structure. We end this section with a short comparison of $\bdcim{n}{2}$, which is the first nontrivial degree-bounded CIM polytope, and $\cim_n$ for small $n$.   

\begin{example}
In the case that $n  = 4$, $\cim_4$ has $2^{n}-n-1 = 11$ coordinates while $\bdcim{4}{2}$ has 10, since $\bc_\cg(1234) = 0$ for all $\cg \in \dags{4}{2}$. Thus, after forgetting all $0$ coordinates, both polytopes are full-dimensional in their ambient spaces of dimensions 11 and 10, respectively. However, their f-vectors, that is, the numbers of faces they have in each dimension from $-1$ (the empty face) to the dimension of the polytope, are
\begin{align*}
&\bdcim{4}{2} :& &(1, 150, 2508, 14387, 39913, 61686, 56250, 30463, 9416, 1495, 94, 1) \\ 
&\cim_4 :& &(1, 185, 4259, 31859, 113320, 225274, 269107, 197894, 88479, 22814, 3013, 154, 1).
\end{align*}
Thus, we can see that $\cim_{4,2}$ is already a great deal simpler with only 150 vertices and 94 facet-defining inequalities compared to $\cim_4$ which has 185 vertices and 154 facets. 

For $n = 5$, this difference only becomes more stark. $\cim_{5,2}$ has only 4066 vertices and has dimension $\binom{n}{2} + \binom{n}{3} = 20$ while $\cim_5$ has 8782 vertices and is 26-dimensional. Lastly, for $n = 6$, $\bdcim{6}{2}$ has 185452 vertices and is $35$-dimensional, while $\cim_6$ is 57-dimensional with 1067825 vertices \cite{enumerating_mecs}. 
\end{example}

It is clear from the example above that while $\bdcim{n}{2}$ is still a very large polytope, it is much smaller and more manageable than the full characteristic imset polytope.

\section{A Facet of $\bdcim{n}{2}$ for Each Connected Graph}
\label{sec:path-ineqs}
\noindent In this section, we provide a large family of inequalities that define facets of $\bdcim{n}{2}$. Each facet is naturally associated to an undirected connected graph $G$ on $m \leq n$ vertices, which records which coordinates $x_{ij}$ appear in the facet-defining inequality. We begin with the following notation, which will be used throughout this paper to define inequalities.

If $\bc_\cg$ is a vertex of $\bdcim{n}{2}$, then $\bc_\cg(S) = 0$ for all subsets $S$ with $|S|\notin\{2,3\}$. Thus, we can naturally consider $\bdcim{n}{2} \subset \rr^{\binom{n}{2} + \binom{n}{3}}$. We write these nonzero coordinates as $x_{ij}$ and $x_{ijk}$ and call them the \emph{edge variables} and \emph{collider variables}, respectively. We write $\be_{ij}$ and $\be_{ijk}$ for the corresponding standard basis vectors of $\rr^{\binom n2+\binom n3}$, so that, for example, the characteristic imset of the graph with all isolated nodes except for the single collider $i\to k\leftarrow j$ is $\be_{ik}+\be_{jk}+\be_{ijk}$. For a set $A$ of coordinates we write $\be_A=\sum_{S\in A}\be_S$ for its indicator vector. When $\cg$ is a DAG, we write $\skel(\cg)$ for its undirected skeleton. For an undirected graph $G$ on $m\leq n$ vertices we write $E(G)$ for its edge set.
We let $C(G)$ denote the set of triples of vertices that induce a path with two edges or a triangle in $G$. That is,
$$C(G) = \big\{\{i,j,k\}: \# \big( \{ ij, ik, jk \} \cap E(G) \big) \geq 2 \big\}.$$
Finally, we let $T(G)$ denote the set of vertices that induce triangles in $G$. In other words,
$$T(G) = \{\{i,j,k\}: ij, ik, jk \in E(G)\}.$$

We also recall the concept of the \emph{incidence matrix of a directed graph}, as given, for instance, in \cite[Section 8.3]{godsil2013algebraic}. 

\begin{definition}
    For a directed graph $\cg$, where each edge $i\rightarrow j$ has a fixed orientation, we define the associated \emph{incidence matrix} $D(\cg)$ as follows. $D(\cg)$ is the $\{0,\pm1\}$-matrix with rows indexed by the vertices of $\cg$ and columns by the edges, such that for edge $e$ and vertex $v$, the $ve$ entry is equal to $+1$ if the vertex $v$ is the head of $e$; $-1$ if $v$ is the tail of $e$; and $0$ otherwise.
\end{definition}

The following theorem regarding the rank of incidence matrices will be useful in the remainder of this section.

\begin{theorem}[{\cite[Theorem 8.3.1]{godsil2013algebraic}}]\label{thm:incidence matrix}
Let $\cg$ be a directed graph on $n$ vertices with $c$ connected components. Then, the rank of the incidence matrix $D(\cg)$ is equal to $n-c$.  
\end{theorem}

\begin{example}\label{eg: chordal example}
Direct computation shows that the following inequality defines a facet of $\bdcim{4}{2}$:
    \begin{equation}\label{eq: chordal example}-(x_{12}+x_{13}+x_{23}+x_{24}+x_{34})+2(x_{123}+x_{234}) + (x_{124}+x_{134})\leq 0.\end{equation}
The edge variables in this inequality can naturally be thought of as indexed by the edges of an undirected graph $G$ pictured in Figure \ref{fig: chordal example}. Note that the subscripts of collider variables $x_{123}$ and $x_{234}$ that appear with coefficient $2$ are exactly the triples in $T(G)$ and the variables $x_{124}$ and $x_{134}$ that appear with coefficient $1$ are the triples in $C(G) \setminus T(G)$. This is the general construction that we will present in Theorem \ref{thm:conn graph}.

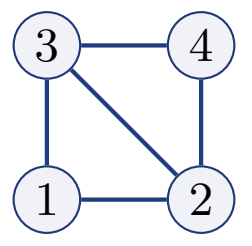
\begin{figure}[ht]
\centering
\setlength{\tabcolsep}{5pt}
\begin{tabular}{c}
\begin{tikzpicture}[>=stealth,x=0.85cm,y=0.85cm,every node/.style={circle,draw=cimnode,fill=cimnodefill,inner sep=0.6pt,minimum size=4.4mm,font=\scriptsize}]
\node (n1) at (0,0) {1};
\node (n2) at (1.2,0) {2};
\node (n3) at (0,1.2) {3};
\node (n4) at (1.2,1.2) {4};
\draw[thick,cimedge] (n1) -- (n2);
\draw[thick,cimedge] (n1) -- (n3);
\draw[thick,cimedge] (n2) -- (n3);
\draw[thick,cimedge] (n2) -- (n4);
\draw[thick,cimedge] (n3) -- (n4);
\end{tikzpicture} \\
 \\[4pt]
\end{tabular}\\
\caption{The undirected graph $G$ corresponding to the facet-defining inequality in Example~\ref{eg: chordal example}.}
\label{fig: chordal example}
\end{figure}

By evaluating the left-hand side of this inequality on all vertices of $\bdcim{4}{2}$, we see that it is indeed valid over the entire polytope. Now let $F$ be the face of $\bdcim{4}{2}$ defined by this inequality. It remains to show that $F$ is a facet. Since $F$ contains the origin, it suffices to find vertices of $F$ whose imsets have rank $\dim(\bdcim{4}{2})-1=9$. We take the zero vector together with the imsets of the $11$ DAGs in Figure~\ref{fig: clique example vertex graphs}: for each triangle $123$ and $234$ of $G$, the three DAGs consisting of one collider on these vertices, for each of the paths $124$ and $134$ the two colliders whose center is not the middle vertex of the path, and the single edge $1\to4$, which is not an edge of $G$. Let $M$ be the matrix whose $12$ columns are the imsets of these DAGs together with the zero vector.

\begin{figure}[ht]
\centering
\setlength{\tabcolsep}{5pt}
\begin{tabular}{cccc}
\begin{tikzpicture}[>=stealth,x=0.85cm,y=0.85cm,every node/.style={circle,draw=cimnode,fill=cimnodefill,inner sep=0.6pt,minimum size=4.4mm,font=\scriptsize}]
\node[coll] (n1) at (0,0) {1};
\node (n2) at (1.2,0) {2};
\node (n3) at (0,1.2) {3};
\node (n4) at (1.2,1.2) {4};
\draw[->,thick,cimedge] (n2) -- (n1);
\draw[->,thick,cimedge] (n3) -- (n1);
\end{tikzpicture} &
\begin{tikzpicture}[>=stealth,x=0.85cm,y=0.85cm,every node/.style={circle,draw=cimnode,fill=cimnodefill,inner sep=0.6pt,minimum size=4.4mm,font=\scriptsize}]
\node (n1) at (0,0) {1};
\node[coll] (n2) at (1.2,0) {2};
\node (n3) at (0,1.2) {3};
\node (n4) at (1.2,1.2) {4};
\draw[->,thick,cimedge] (n1) -- (n2);
\draw[->,thick,cimedge] (n3) -- (n2);
\end{tikzpicture} &
\begin{tikzpicture}[>=stealth,x=0.85cm,y=0.85cm,every node/.style={circle,draw=cimnode,fill=cimnodefill,inner sep=0.6pt,minimum size=4.4mm,font=\scriptsize}]
\node (n1) at (0,0) {1};
\node (n2) at (1.2,0) {2};
\node[coll] (n3) at (0,1.2) {3};
\node (n4) at (1.2,1.2) {4};
\draw[->,thick,cimedge] (n1) -- (n3);
\draw[->,thick,cimedge] (n2) -- (n3);
\end{tikzpicture} &
\begin{tikzpicture}[>=stealth,x=0.85cm,y=0.85cm,every node/.style={circle,draw=cimnode,fill=cimnodefill,inner sep=0.6pt,minimum size=4.4mm,font=\scriptsize}]
\node (n1) at (0,0) {1};
\node[coll] (n2) at (1.2,0) {2};
\node (n3) at (0,1.2) {3};
\node (n4) at (1.2,1.2) {4};
\draw[->,thick,cimedge] (n3) -- (n2);
\draw[->,thick,cimedge] (n4) -- (n2);
\end{tikzpicture} \\
{\scriptsize $\cg_{123}^{(1)}$} & {\scriptsize $\cg_{123}^{(2)}$} & {\scriptsize $\cg_{123}^{(3)}$} & {\scriptsize $\cg_{234}^{(2)}$} \\[4pt]
\end{tabular}\\
\begin{tabular}{cccc}
\begin{tikzpicture}[>=stealth,x=0.85cm,y=0.85cm,every node/.style={circle,draw=cimnode,fill=cimnodefill,inner sep=0.6pt,minimum size=4.4mm,font=\scriptsize}]
\node (n1) at (0,0) {1};
\node (n2) at (1.2,0) {2};
\node[coll] (n3) at (0,1.2) {3};
\node (n4) at (1.2,1.2) {4};
\draw[->,thick,cimedge] (n2) -- (n3);
\draw[->,thick,cimedge] (n4) -- (n3);
\end{tikzpicture} &
\begin{tikzpicture}[>=stealth,x=0.85cm,y=0.85cm,every node/.style={circle,draw=cimnode,fill=cimnodefill,inner sep=0.6pt,minimum size=4.4mm,font=\scriptsize}]
\node (n1) at (0,0) {1};
\node (n2) at (1.2,0) {2};
\node (n3) at (0,1.2) {3};
\node[coll] (n4) at (1.2,1.2) {4};
\draw[->,thick,cimedge] (n2) -- (n4);
\draw[->,thick,cimedge] (n3) -- (n4);
\end{tikzpicture} &
\begin{tikzpicture}[>=stealth,x=0.85cm,y=0.85cm,every node/.style={circle,draw=cimnode,fill=cimnodefill,inner sep=0.6pt,minimum size=4.4mm,font=\scriptsize}]
\node[coll] (n1) at (0,0) {1};
\node (n2) at (1.2,0) {2};
\node (n3) at (0,1.2) {3};
\node (n4) at (1.2,1.2) {4};
\draw[->,thick,cimedge] (n2) -- (n1);
\draw[->,thick,cimedge] (n4) -- (n1);
\end{tikzpicture} &
\begin{tikzpicture}[>=stealth,x=0.85cm,y=0.85cm,every node/.style={circle,draw=cimnode,fill=cimnodefill,inner sep=0.6pt,minimum size=4.4mm,font=\scriptsize}]
\node (n1) at (0,0) {1};
\node (n2) at (1.2,0) {2};
\node (n3) at (0,1.2) {3};
\node[coll] (n4) at (1.2,1.2) {4};
\draw[->,thick,cimedge] (n1) -- (n4);
\draw[->,thick,cimedge] (n2) -- (n4);
\end{tikzpicture} \\
{\scriptsize $\cg_{234}^{(3)}$} & {\scriptsize $\cg_{234}^{(4)}$} & {\scriptsize $\cg_{124}^{(1)}$} & {\scriptsize $\cg_{124}^{(4)}$} \\[4pt]
\end{tabular}\\
\begin{tabular}{cccc}
\begin{tikzpicture}[>=stealth,x=0.85cm,y=0.85cm,every node/.style={circle,draw=cimnode,fill=cimnodefill,inner sep=0.6pt,minimum size=4.4mm,font=\scriptsize}]
\node[coll] (n1) at (0,0) {1};
\node (n2) at (1.2,0) {2};
\node (n3) at (0,1.2) {3};
\node (n4) at (1.2,1.2) {4};
\draw[->,thick,cimedge] (n3) -- (n1);
\draw[->,thick,cimedge] (n4) -- (n1);
\end{tikzpicture} &
\begin{tikzpicture}[>=stealth,x=0.85cm,y=0.85cm,every node/.style={circle,draw=cimnode,fill=cimnodefill,inner sep=0.6pt,minimum size=4.4mm,font=\scriptsize}]
\node (n1) at (0,0) {1};
\node (n2) at (1.2,0) {2};
\node (n3) at (0,1.2) {3};
\node[coll] (n4) at (1.2,1.2) {4};
\draw[->,thick,cimedge] (n1) -- (n4);
\draw[->,thick,cimedge] (n3) -- (n4);
\end{tikzpicture} &
\begin{tikzpicture}[>=stealth,x=0.85cm,y=0.85cm,every node/.style={circle,draw=cimnode,fill=cimnodefill,inner sep=0.6pt,minimum size=4.4mm,font=\scriptsize}]
\node (n1) at (0,0) {1};
\node (n2) at (1.2,0) {2};
\node (n3) at (0,1.2) {3};
\node (n4) at (1.2,1.2) {4};
\draw[->,thick,cimedge] (n1) -- (n4);
\end{tikzpicture} &
 \\
{\scriptsize $\cg_{134}^{(1)}$} & {\scriptsize $\cg_{134}^{(4)}$} & {\scriptsize $\cg_{14}$} &  \\[4pt]
\end{tabular}\\
\caption{The $11$ DAGs whose imsets are the nonzero columns of $M$ in Example~\ref{eg: chordal example}. Purple vertices are colliders.}
\label{fig: clique example vertex graphs}
\end{figure}
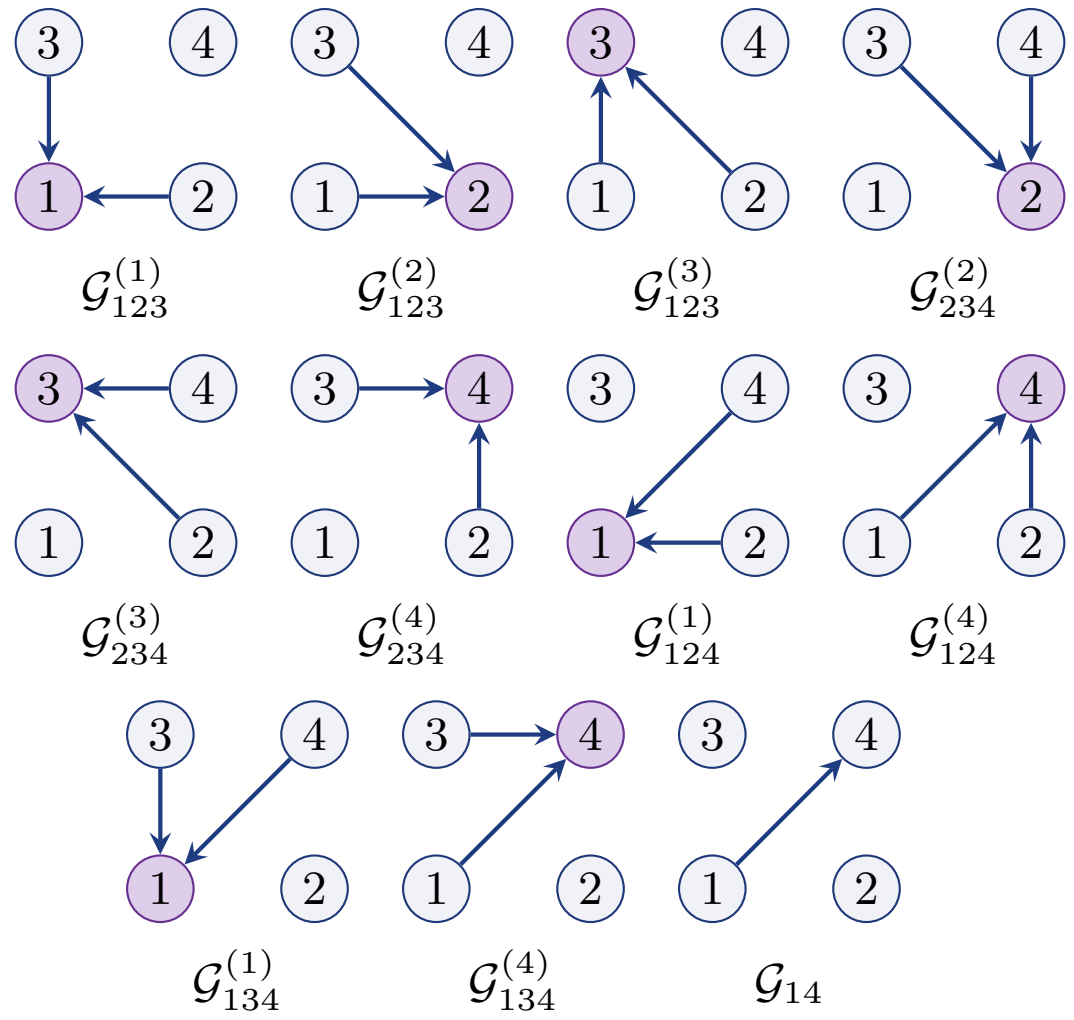

Subtract the column $\be_{14}$ from the four columns supported on the row $14$, and, for each of the columns supported on a given triple $ijk$, subtract the first of those columns from the others. The resulting columns are of two kinds. Five of them - the first column corresponding to each triple, and $\be_{14}$ - each are supported on that triple (or the $14$ row in the case of $\be_{14}$), and these rows are zero in all other columns. The remaining six columns are supported on the edge rows $12,13,23,24,34$ only. They form the matrix:
\[
D=\begin{array}{c|cccccc}
   & 12\!\to\!23 & 13\!\to\!23 & 23\!\to\!34 & 24\!\to\!34 & 12\!\to\!24 & 13\!\to\!34\\ \hline
12 & -1 & 0 & 0 & 0 & -1 & 0\\
13 & 0 & -1 & 0 & 0 & 0 & -1\\
23 & 1 & 1 & -1 & 0 & 0 & 0\\
24 & 0 & 0 & 0 & -1 & 1 & 0\\
34 & 0 & 0 & 1 & 1 & 0 & 1
\end{array}.
\]
This is the incidence matrix of the DAG $\ch$ in Figure~\ref{fig: incidence digraph}. The columns are labeled by the corresponding edges of $\ch$. Since $\ch$ is connected on $5$ vertices, Theorem~\ref{thm:incidence matrix} gives $\rank(D)=4$, so $M$ has rank $5+4=9$.

\begin{figure}[ht]
\centering
\setlength{\tabcolsep}{5pt}
\begin{tabular}{c}
\begin{tikzpicture}[>=stealth,x=0.8cm,y=0.8cm,every node/.style={circle,draw=cimnode,fill=cimnodefill,inner sep=0.6pt,minimum size=6.5mm,font=\scriptsize}]
\node (n12) at (0,3.0) {12};
\node (n13) at (0,0.2) {13};
\node (n23) at (2,1.6) {23};
\node (n24) at (2,3.0) {24};
\node (n34) at (4,1.6) {34};
\draw[->,thick,cimedge] (n12) -- (n23);
\draw[->,thick,cimedge] (n13) -- (n23);
\draw[->,thick,cimedge] (n23) -- (n34);
\draw[->,thick,cimedge] (n24) -- (n34);
\draw[->,thick,cimedge] (n12) -- (n24);
\draw[->,thick,cimedge] (n13) -- (n34);
\end{tikzpicture} \\
 \\[4pt]
\end{tabular}\\
\caption{The directed graph $\ch$ whose incidence matrix is the matrix above; its vertices are the edges of $G$.}
\label{fig: incidence digraph}
\end{figure}
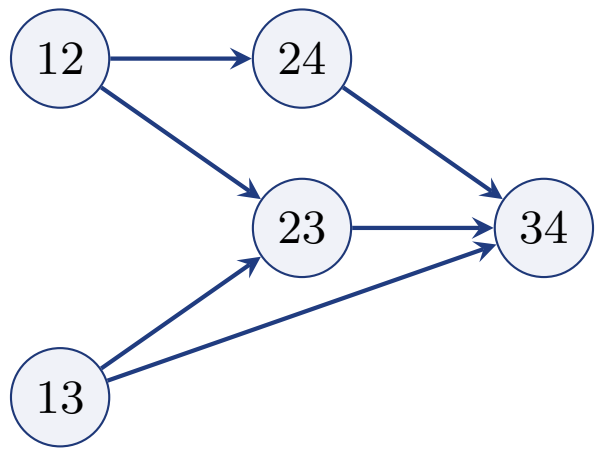
\end{example}

This example illustrates one facet of $\cim_{4,2}$ which is naturally indexed by a connected graph on 4 vertices. The following lemma shows that this inequality also defines a facet of $\cim_{n,2}$ for all $n \geq 4$.

\begin{lemma}\label{lem: building facets}
    Let $\ba\in\rr^{\binom n2+\binom n3}$ and $a_0\in\rr$ be such that $\ba\cdot\bx\leq a_0$ is a valid inequality on $\bdcim{n}{2}$, and let $\ba'\in\rr^{\binom{n+1}2+\binom{n+1}3}$ be the vector obtained from $\ba$ by setting $\ba'_S=0$ for every coordinate $S$ that contains $n+1$. If $F = \bdcim{n}{2} \cap \{\bx\mid\ba\cdot\bx= a_0\}$ is a facet of $\bdcim{n}{2}$, then $F' = \bdcim{n+1}{2}\cap\{\bx\mid\ba'\cdot\bx= a_0\}$ is a facet of $\bdcim{n+1}{2}$.
\end{lemma}
\begin{proof}
    For $\cg\in\dags{n+1}{2}$, the restriction of $\bc_\cg$ to the coordinates indexed by subsets of $[n]$ is the imset of the induced subgraph $\cg[[n]]\in\dags{n}{2}$. Since $\ba'$ vanishes on all other coordinates, $\ba'\cdot\bc_\cg=\ba\cdot\bc_{\cg[[n]]}\leq a_0$. Hence $\ba'\cdot\bx\leq a_0$ is valid on $\bdcim{n+1}{2}$ and defines a face $F'$, and it remains to show that $F'$ is a facet.

    Let $M$ be the defining matrix of $F$. Since $F$ is a facet and $\bdcim{n}{2}$ is full-dimensional, $\rank(\widehat{M})=\dim(\bdcim{n}{2})=\binom n2+\binom n3$ where $\widehat{M}$ is as defined in Definition \ref{def: hat}.. Fix a DAG $\cg_0\in\dags{n}{2}$ with $\bc_{\cg_0}\in F$.  We regard DAGs on $[n]$ as DAGs on $[n+1]$ in which $n+1$ is isolated; their imsets are obtained by adding zeros in the coordinates that contain $n+1$. Let $M'$ be the matrix whose columns are the imsets of the following DAGs in $\dags{n+1}{2}$.
    \begin{enumerate}
        \item For each column $\bc_\cg$ of $M$, the DAG $\cg$ with $n+1$ added as an isolated vertex.
        \item For each $i\in[n]$, the DAG $\cg_0+(i\to n+1)$, with imset $\bc_{\cg_0}+\be_{i\,n+1}$.
        \item For each pair $i,j\in[n]$, the DAG $\cg_0+(i\to n+1)+(j\to n+1)$, with imset $\bc_{\cg_0}+\be_{i\,n+1}+\be_{j\,n+1}+\be_{ij\,n+1}$.
    \end{enumerate}
    In each of these DAGs, the vertex $n+1$ is a sink with at most two parents, so they all lie in $\dags{n+1}{2}$. Their restrictions to the coordinates indexed by subsets of $[n]$ are $\bc_\cg$ or $\bc_{\cg_0}$, which lie on $F$, so all of them satisfy $\ba'\cdot\bx=a_0$ and lie on $F'$.

    Order the rows of $\widehat{M'}$ so that the row of ones and the coordinates indexed by subsets of $[n]$ come first, followed by the coordinates $\{i,n+1\}$ and then the coordinates $\{i,j,n+1\}$. Order the columns as in (1), (2), (3). Then
    \[
    \widehat{M'}=\begin{pmatrix}\widehat{M}&*\\ 0&N\end{pmatrix},
    \]
    where $N$ is the square matrix of size $n+\binom n2$ formed by the coordinates containing $n+1$ of the columns in (2) and (3). The column of (2) for $i$ has a single $1$ in row $\{i,n+1\}$. The column of (3) for $\{i,j\}$ has a $1$ in row $\{i,j,n+1\}$, and otherwise only entries in rows $\{i,n+1\}$ and $\{j,n+1\}$, which come earlier. So $N$ is upper triangular with ones on the diagonal. Since the lower left block of $\widehat{M'}$ is zero,
    \[
    \rank(\widehat{M'})\geq\rank(\widehat{M})+\rank(N)=\tbinom n2+\tbinom n3+n+\tbinom n2=\tbinom{n+1}2+\tbinom{n+1}3=\dim(\bdcim{n+1}{2}).
    \]
    Therefore $F'$ is a facet of $\bdcim{n+1}{2}$.
\end{proof}

To prove the main theorem of this section, we require the following notation. Let $\cg$ be a DAG on the same vertex set as $G$. Denote by $\cg [G]$ the subgraph of $\cg$ obtained by restricting the skeleton of $\cg$ to the edges in $E(G)$. Let $O(\cg)$ denote the set of all edges of $\cg$ which do not belong to any collider. We say that an \emph{unordered} pair $ij$ belongs to $O(\cg)$ if either $i \rightarrow j$ or $j \rightarrow i$ belongs to $O(\cg)$.

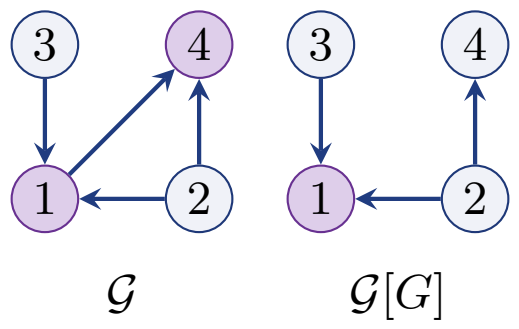
\begin{figure}[ht]
\centering
\setlength{\tabcolsep}{5pt}
\begin{tabular}{cc}
\begin{tikzpicture}[>=stealth,x=0.85cm,y=0.85cm,every node/.style={circle,draw=cimnode,fill=cimnodefill,inner sep=0.6pt,minimum size=4.4mm,font=\scriptsize}]
\node[coll] (n1) at (0,0) {1};
\node (n2) at (1.2,0) {2};
\node (n3) at (0,1.2) {3};
\node[coll] (n4) at (1.2,1.2) {4};
\draw[->,thick,cimedge] (n2) -- (n1);
\draw[->,thick,cimedge] (n2) -- (n4);
\draw[->,thick,cimedge] (n3) -- (n1);
\draw[->,thick,cimedge] (n1) -- (n4);
\end{tikzpicture} &
\begin{tikzpicture}[>=stealth,x=0.85cm,y=0.85cm,every node/.style={circle,draw=cimnode,fill=cimnodefill,inner sep=0.6pt,minimum size=4.4mm,font=\scriptsize}]
\node[coll] (n1) at (0,0) {1};
\node (n2) at (1.2,0) {2};
\node (n3) at (0,1.2) {3};
\node (n4) at (1.2,1.2) {4};
\draw[->,thick,cimedge] (n2) -- (n1);
\draw[->,thick,cimedge] (n2) -- (n4);
\draw[->,thick,cimedge] (n3) -- (n1);
\end{tikzpicture} \\
{\scriptsize $\cg$} & {\scriptsize $\cg[G]$} \\[4pt]
\end{tabular}\\
\caption{The DAG $\cg$ from Example~\ref{eg: restriction notation} and its restriction $\cg[G]$. Purple vertices are colliders.}
\label{fig: restriction notation}
\end{figure}

\begin{example}\label{eg: restriction notation}
    Let $G$ be as in Example \ref{eg: chordal example} and let $\cg$ be the DAG on the left of Figure \ref{fig: restriction notation}. Then $\cg[G]$ is the graph obtained from $\cg$ by deleting all edges of the skeleton of $\cg$ that are not in $G$; since $14$ is the only such edge, we obtain the DAG on the right of Figure \ref{fig: restriction notation}. Since all edges of $\cg$ belong to a collider, we have $O(\cg) = \emptyset$. However, since $1 \rightarrow 4$ is not an edge of $\cg[G]$, we have $O(\cg[G]) = \{2 \rightarrow 4 \}$.
\end{example}

\begin{theorem}\label{thm:conn graph}
     Let $G$ be a connected undirected graph on $[m]$.  Then, the following inequality defines a facet of $\bdcim{n}{2}$, for $n\geq m$. 
     \begin{equation}\label{eq: every undirected graph inequality} -\left(\sum_{(i,j)\in E(G)}x_{ij}\right) + 2\left(\sum_{(i,j,k)\in T(G)}x_{ijk}\right) + \sum_{(i,j,k)\in C(G)\setminus T(G)}x_{ijk} \leq 0 .\end{equation}
We call \eqref{eq: every undirected graph inequality} the \emph{path inequality} of $G$.
\end{theorem}
\begin{proof}
         By Lemma~\ref{lem: building facets}, applied $n-m$ times, it suffices to show that the inequality \eqref{eq: every undirected graph inequality} is valid on $\bdcim{m}{2}$ and defines a facet of $\bdcim{m}{2}$; extending the functional by zeros then gives a facet of $\bdcim{n}{2}$ for every $n\geq m$. Throughout the proof we therefore work with DAGs on the vertex set $[m]$ of $G$.

        \textbf{Validity:} 
        We check that the inequality (\ref{eq: every undirected graph inequality}) is valid on $\bdcim{m}{2}$ by evaluating the left-hand side on all vertices of the polytope. 
   Let $\bc_\cg$ be the characteristic imset of a DAG $\cg\in\dags{m}{2}$. We consider the value of 
   \[
   -\left(\sum_{(i,j)\in E(G)}\bc_\cg (ij)\right) + 2\left(\sum_{(i,j,k)\in T(G)}\bc_\cg (ijk)\right) + \sum_{(i,j,k)\in C(G)\setminus T(G)}\bc_\cg(ijk).
   \]
   Since $\cg\in\dags{m}{2}$, every node of $\cg$ has at most two parents. Hence, each edge $i \rightarrow j$ is in at most one collider.  
   
   If the entry  $\bc_\cg(ijk) = 1$, then without loss of generality $i \rightarrow j \leftarrow k$ forms a collider in $\cg$. Thus $\bc_\cg(ij) = \bc_\cg(jk) = 1$. If $ijk$ forms a triangle in $G$, then $-\bc_\cg(ij)$ and $-\bc_\cg(jk)$ both appear in the inequality and we have $-\bc_\cg(ij) - \bc_\cg(jk) + 2\bc_\cg(ijk) = 0$. If $ijk\in C(G)\setminus T(G)$, then at least two of the pairs $ij,ik,jk$ are edges of $G$, so at least one of $-\bc_\cg(ij)$ and $-\bc_\cg(jk)$ appears in the inequality. If both do, we have $-\bc_\cg(ij) - \bc_\cg(jk) + \bc_\cg(ijk) = -1$. If only one does, then suppose without loss of generality that it is $ij$; then we have $-\bc_\cg(ij) + \bc_\cg(ijk) = 0$. If $ijk\notin C(G)$, then $x_{ijk}$ does not appear in the inequality. Since each edge of $\cg$ belongs to at most one collider, we can obtain the value of the left-hand side of (\ref{eq: every undirected graph inequality}) evaluated at $\bc_\cg$ as
   \begin{eqnarray*}
       \sum_{\substack{i \rightarrow j \leftarrow k \in \mathrm{coll}(\cg) \\ ijk \in T(G)}} \big(-\bc_\cg(ij) - \bc_\cg(jk) + 2\bc_\cg(ijk) \big) + \sum_{\substack{i \rightarrow j \leftarrow k \in \mathrm{coll}(\cg) \\ ijk \in C(G) \setminus T(G) \\ ij, jk \in E(G)}} \big(-\bc_\cg(ij) - \bc_\cg(jk) + \bc_\cg(ijk) \big)  \\
       + \sum_{\substack{i \rightarrow j \leftarrow k \in \mathrm{coll}(\cg) \\ ijk \in C(G) \setminus T(G) \\ \#(\{ij, jk\} \cap E(G)) = 1}} \big( -\bc_\cg(ij) + \bc_\cg(ijk) \big) - \sum_{ij \in R} \bc_\cg(ij),
   \end{eqnarray*}
   where, in the third sum, $ij$ denotes the edge of the collider that lies in $E(G)$, and $R$ is the set of edges of $\cg[G]$ that do not lie in a collider $i\to j\leftarrow k$ of $\cg$ with $ijk\in C(G)$. Since each summand is less than or equal to $0$, by the above discussion, the value of this linear functional at $\bc_\cg$ is also less than or equal to $0$ as needed.

        \textbf{Facet-defining:} 
        Let $F$ be the face of $\bdcim{m}{2}$ defined by the inequality \eqref{eq: every undirected graph inequality}.
        We choose a subset of the vertices of $F$ to be the columns of a matrix $M$ and prove that the affine hull of these columns has dimension one less than that of $\bdcim{m}{2}$. Let $M$ be the matrix whose columns are the following characteristic imsets:

         \begin{enumerate}
        \item For each triple $ijk\in T(G)$ with $i < j <k$, we have three columns. The first is the imset of $\cg_{ijk}^{(i)} = j\rightarrow i \leftarrow k$; the second of $\cg_{ijk}^{(j)} = i\rightarrow j \leftarrow k$; and the third of $\cg_{ijk}^{(k)} = i\rightarrow k \leftarrow j$. In each of these graphs, all vertices in $[m] \setminus \{i,j,k\}$ are isolated.
        \item  
        For each triple $ijk \in C(G) \setminus T(G)$, where $ij,jk\in E(G)$, $ik\notin E(G)$ and $i < k$, we obtain two more columns. These are the imsets of $\cg_{ijk}^{(i)} = j\rightarrow i \leftarrow k$ and $\cg_{ijk}^{(k)} = i\rightarrow k \leftarrow j$. Again, each of these is a graph on the vertex set $[m]$ with all other nodes isolated.
        \item For each pair $ij \notin E(G)$, the imset $\be_{ij}$ of the DAG with the single edge $i\to j$. 
        \item For each triple $ijk\notin C(G)$, at most one of its pairs is an edge of $G$, so we may label it so that $ik,jk\notin E(G)$; we take the imset $\be_{ik}+\be_{jk}+\be_{ijk}$ of the collider $i\to k\leftarrow j$.
    \end{enumerate}

    Since $F$ contains the origin, the dimension of $F$ is at least $\rank(M)$. In fact, we will show that $M$ has rank $\dim(\bdcim{m}{2}) - 1$, which shows that $F$ is a facet. We perform column operations on $M$ in the following ways to make its rank easier to compute. 
    For each imset in item (2), associated to the graphs $\cg^{(i)}_{ijk}$ and $\cg^{(k)}_{ijk}$ where $ik$ is not an edge of $G$, 
    we subtract the column $\be_{ik}$ to obtain $\be_{ij}+\be_{ijk}$ and $\be_{jk}+\be_{ijk}$. Similarly, for each triple $ijk\notin C(G)$ in item (4), we subtract the columns $\be_{ik}$ and $\be_{jk}$ to obtain $\be_{ijk}$. These operations do not change the rank of the matrix. Let $M'$ be the matrix obtained from $M$ by these operations.

    We arrange the rows and columns of $M'$ so that all rows corresponding to $ij \not\in E(G)$ and $ijk \not\in C(G)$ are at the bottom of $M'$ and all corresponding indicator vector columns $\be_{ij}$ and $\be_{ijk}$ are on the right in the same order. 
    Thus  $M'$ is block diagonal with lower block equal to the identity of rank $\dim(\bdcim{m}{2})-|E(G)|-|C(G)|$. The upper block $U$ consists of the remaining columns.
    
    For each triple $ijk \in T(G)$, we have three columns in $U$ whose restrictions to the $ij, ik, jk$ and $ijk$ rows are
   \[ \begin{blockarray}{cccc}
     & \cg_{ijk}^{(i)} & \cg_{ijk}^{(j)} & \cg_{ijk}^{(k)}\\
        \begin{block}{c(ccc)}
           ij & 1 & 1 & 0 \\ 
           ik& 1 & 0 &1 \\ 
           jk&0 & 1 & 1\\ 
           ijk &1 & 1 & 1\\
        \end{block}
    \end{blockarray}. \]
    We can then perform simple row and column operations as follows: 
     \begin{align*}
       \xrightarrow[\cg_{ijk}^{(k)}-\cg_{ijk}^{(i)} \text{ on } \cg_{ijk}^{(k)}]{\cg_{ijk}^{(j)}-\cg_{ijk}^{(i)} \text{ on } \cg_{ijk}^{(j)}} 
        \begin{pmatrix}
            1 & 0 & -1 \\
            1 & -1 & 0 \\
            0 & 1 & 1 \\
            1 & 0 & 0
        \end{pmatrix} \xrightarrow[r_{ik}-r_{ijk} \text{ on } r_{ik}]{r_{ij}-r_{ijk} \text{ on } r_{ij}} 
        \begin{pmatrix}
            0 & 0 & -1 \\
            0 & -1 & 0 \\
            0 & 1 & 1 \\
            1 & 0 & 0
        \end{pmatrix}.
    \end{align*}
Note that these operations do not affect any other entries of this matrix, since the $ijk$ row and the $\cg_{ijk}^{(i)}, \cg_{ijk}^{(j)}$ and $\cg_{ijk}^{(k)}$ columns have all zeros outside of these entries. 
    For the triples $ijk \in C(G) \setminus T(G)$ where $ij, jk \in E(G)$ and $ik \not\in E(G)$, we have two columns of $U$ which are supported only on the $ij, jk$ and $ijk$ rows. Restricting to these rows, these columns are
    \[\begin{blockarray}{ccc}
     & \cg_{ijk}^{(i)} & \cg_{ijk}^{(k)}\\
        \begin{block}{c(cc)}
            ij & 1 & 0\\
            jk & 0 & 1\\
            ijk & 1 & 1\\
        \end{block}
    \end{blockarray}\]
     We can subtract the $\cg_{ijk}^{(i)}$ column from the $\cg_{ijk}^{(k)}$ column and the $ijk$ row from the $ij$ row to obtain $\begin{pmatrix}
        0& -1 \\ 0 & 1 \\ 1 & 0
    \end{pmatrix}$.
    For the same reason as the previous case, these row and column operations do not affect any other entries of the matrix.

    We can thus think of $U$ to be formed of blocks of columns, each of which is either three or two columns wide. Also, note that we can change the ordering of the columns and rows to ensure that within each such block, the only edge variable row that has a positive entry in this block is the lowest edge variable row of the block. 
    
    Thus, we can reorder the rows and columns of $U$ to be a block diagonal matrix, with lower block $I_{|C(G)|}$. The upper block can be realized as the incidence matrix of a directed graph $\cg'$ on vertices indexed by $ij$, where each $ij$ is an edge of $G$. The edges of $\cg'$ are as follows. For each triple $ijk \in T(G)$, we have two edges: $ij\rightarrow jk$ and $ik\rightarrow jk$. For each $ijk \in C(G)\setminus T(G)$, where $ik\notin E(G)$, there is one edge $ij \rightarrow jk$. Any two edges $ij$ and $jk$ of $G$ that share a vertex lie in a common triple $ijk\in C(G)$, and so they are joined by a path in $\cg'$, either directly by an edge or through $ik$. Since $G$ is connected, so is its line graph, and hence $\cg'$ is connected. Thus, by Theorem~\ref{thm:incidence matrix}, the rank of the upper block, which is also the incidence matrix of $\cg'$, is $|E(G)|-1$.

    Therefore, $F$ is a facet of $\bdcim{m}{2}$, as

    \begin{align*}
        \dim(F)\geq\rank(M') & = \dim(\bdcim{m}{2})-|E(G)|-|C(G)| + \rank(U)\\
        & = \dim(\bdcim{m}{2})-|E(G)|-|C(G)| + |C(G)| + |E(G)| - 1\\
        & = \dim(\bdcim{m}{2}) -1. 
    \end{align*}
    \end{proof}

\section{Additional Facets from Trees} 
\label{sec:star-tree-ineqs}

We now describe another family of facets that arise from trees. Let $T$ be a tree on $[m]$, $m\geq2$, with no vertex of degree $2$. The vertices of $T$ of degree $1$ are called \emph{leaves} and all others are \emph{internal}, namely those of degree $\geq3$. Write $N(v)=N_T(v)$ for the set of neighbors of $v$ in $T$.

\begin{definition}\label{def:outer}
Let $\Nb=\Nb(T)=\bigcup_{v\text{ internal}}\binom{N(v)}{3}$. The \emph{outer inequality} of $T$ is:
\begin{equation}\label{eq:tree inequality}
\ba_T\cdot \bx:=\sum_{ij\in E(T)}x_{ij}+\sum_{ijk\in\Nb}x_{ijk}\le|E(T)|=m-1 .
\end{equation}
\end{definition}

\noindent The primary goal of this section is to prove the following theorem. 

\begin{theorem}\label{thm:tree-facet}
For every such $T$ be a tree on $m$ vertices with no vertices of degree $2$. Then for every $n\ge m$, the outer inequality \eqref{eq:tree inequality} is valid and facet-defining for $\bdcim{n}{2}$.
\end{theorem}

 Our proof approach is similar to that in the previous section in that we first prove that the outer inequality is valid on $\bdcim{n}{2}$ and then exhibit a specific collection of DAGs whose characteristic imsets lie on the facet and whose affine span has the correct dimension. We first introduce some notation which will be helpful throughout the proof. 

Let $S=E(T)\cup\Nb$ be the support of $\ba_T$; $d=\dim\bdcim m2=\binom m2+\binom m3$; let $\Lambda$ be the set of the other $d-|S|$ coordinates (pairs that are not edges of $T$ and triples not in $\Nb$).
We write $\mathbf 1_T=\sum_{ij\in E(T)}\be_{ij}$.
Let $\overline T$ be the graph obtained from $T$ by adding all pairs of vertices that have a common neighbor in $T$.  We call a DAG $\cg\in\dags{m}{2}$ \emph{tight} if $\ba_T\cdot \bc_\cg=m-1$; the face $F$ defined by \eqref{eq:tree inequality} is the convex hull of the imsets of tight DAGs. Since $T$ is a tree, the following three facts hold; we label them for later reference.
\begin{enumerate}[label=(T\arabic*),ref=(T\arabic*),nosep]
\item\label{item:T1} If $x,y\in N(v)$ then $xy\notin E(T)$. Thus a triple belongs to $\binom{N(v)}{3}$ for at most one $v$.
\item\label{item:T2} If $u\in N(v)$, then there is no triple in $\Nb$ containing both $u$ and $v$.
\item\label{item:T3} If $u\in N(v)$ and $x \in N(u)$, then $x \notin N(v)$.
\end{enumerate}

\begin{lemma}\label{lem: paths in T bar}
Let $x$ be a vertex of $T$ and $u,v\in N(x)$. Every path from $u$ to $v$ in $\overline T$ avoiding $x$ has all its vertices in $N(x)$.
\end{lemma}
\begin{proof}
For $y\in N(x)$, let $T_y$ be the component of $T-x$ containing $y$ and $S_y=V(T_y)\setminus\{y\}$; these $V(T_y)$ partition $V(T)\setminus\{x\}$ and $y$ is the only vertex of $T_y$ adjacent to $x$ in $T$. 

We claim that if $w\in S_y$ and $z\ne x$ is adjacent to $w$ in $\overline T$ then $z\in V(T_y)$. Indeed, either $wz\in E(T)$, or $w,z$ have a common neighbor $t$; if $t\neq x$ the edges $wt,tz$ lie in $T-x$, and if $t=x$ then $w\in N(x)\cap V(T_y)=\{y\}$, contradicting $w\in S_y$. This proves the claim.

So the only vertex outside $S_y$ (other than $x$) adjacent to an element of $S_y$ is $y$. If a $uv$-path $P$ avoiding $x$ had a vertex outside $N(x)$, it would lie in some $S_y$; as $u,v\notin S_y$, $P$ enters and leaves $S_y$, both times through $y$, so it visits $y$ twice, which contradicts that $P$ is a path. This proves the lemma.
\end{proof}

\subsection{Validity}\label{ssec:validity}
Let $\cg\in\dags{n}{2}$. Each vertex $w$ with two parents is a collider corresponding to the triple $\{w\}\cup\mathrm{pa}(w)$. Call $w$ a \emph{support collider} if this triple is in $\Nb$, and let $\mathrm{Col}_T(\cg)$ be the set of support colliders in $\cg$. Then, 
\begin{equation}\label{eq:value}
\ba_T\cdot \bc_\cg=|E\left( \cg[T] \right)|+|\mathrm{Col}_T(\cg)|.
\end{equation}
This depends only on the induced subgraph $\cg\left[ \{1,\dots,m\} \right]$, so by Lemma \ref{lem: building facets}, it suffices to consider the $n=m$ case.

\begin{proposition}\label{prop:step}
If $\mathrm{Col}_T(\cg)\ne\emptyset$, then there is a DAG $\cg'\in\dags{m}{2}$ whose skeleton is contained in $\overline{T}$ with $\ba_T\cdot \bc_{\cg'}=\ba_T\cdot \bc_\cg$ and $|\mathrm{Col}_T(\cg')|<|\mathrm{Col}_T(\cg)|.$
\end{proposition}

In the proof of the proposition, we will make use of the following remark.

\begin{remark}\label{rem:Acyclicity}
   Let $\ch$ be a DAG and suppose that $\ch$ contains a directed path $a \rightsquigarrow b$ for two vertices $a$ and $b$ . Let $\ch'$ be obtained from any subgraph of $\ch$ by adding the edge $a\rightarrow b$. Then $\ch'$ is again a DAG, since otherwise, replacing the added edge of a directed cycle in $\ch'$ by the corresponding path in $\ch$ would give a closed directed walk in $\ch$, which contradicts that $\ch$ is acyclic.
\end{remark}

\begin{proof}[Proof of Proposition~\ref{prop:step}]
Let $\cg \in \dags m2$. Since deleting edges of $\cg$ that are not (undirected) edges of $\overline{T}$ does not change $\ba_T \cdot \bc_\cg$, we may assume that the edges of $\cg$ are contained in the edges of $\overline{T}$.
Let $u_2$ be a support collider of $\cg$ with parents $u_1,u_3$ and $\{u_1,u_2,u_3\}\subseteq N(v_1)$. By~\ref{item:T1}, $u_1u_2,u_3u_2\notin E(T)$. Define $\mathcal{B}=\cg-\{u_1 \rightarrow u_2,u_3 \rightarrow u_2\}$, which has one fewer support collider than $\cg$, in particular the one at $u_2$, and the same number of edges of $T$. 
Below, we construct a digraph $\cg'$ from $\mathcal{B}$ depending on the relationship between $u_2$ and $v_1$ in $\cg$ and show that in each case, $\cg'\in\dags{m}{2}$ and that $\ba_T \cdot \bc_\cg = \ba_T \cdot \bc_{\cg'}$.

\emph{Case 1:} Suppose $u_2$ is a descendant of $v_1$. Then, $v_1 \rightarrow u_2\notin E(\cg)$ since $u_2$ already has two parents $u_1,u_3$ in $\cg$. Let $\cg':= \mathcal{B}+(v_1\rightarrow u_2)$. This is acyclic by Remark~\ref{rem:Acyclicity}, since $\cg$ has a path $v_1\rightsquigarrow u_2$, and $u_2$ has in-degree $1$ in $\cg'$. The DAG $\cg'$ is obtained from $\cg$ by adding the edge $v_1u_2$ of $T$ and deleting a support collider, so $\ba_T\cdot \bc_{\cg'}=\ba_T\cdot \bc_\cg$.

\emph{Case 2:} Suppose $u_2$ is neither a descendant nor an ancestor of $v_1$. We again see that $\cg'=\mathcal{B}+(v_1\to u_2)$ is acyclic as there is no path $u_2\rightsquigarrow v_1$, and the rest of the argument is analogous to Case 1.

\emph{Case 3:} Suppose that $u_2$ is an ancestor of $v_1$. Then, so are $u_1,u_3$. Thus none of $u_1,u_2,u_3$ is a child of $v_1$, and $v_1$ has at least one parent. 

\emph{3a:} First we assume that $|\mathrm{pa}(v_1)|=2$ and $\mathrm{pa}(v_1)=\{x_1,x_2\}\subseteq N(w)$ for some $w\in N(v_1)$. Then, $v_1$ is a support collider. By~\ref{item:T3}, $x_1, x_2\notin N(v_1)$, so $x_1, x_2 \notin\{u_1,u_2,u_3\}$, and $u_1 \rightarrow v_1$ and $u_3\rightarrow v_1$ are not edges of $\cg$. Furthermore, since $T$ has no triangles, $x_1v_1, x_2v_1\notin E(T)$. Set $\cg'=\mathcal{B} + \{ u_1\!\to\! v_1,u_3\!\to\! v_1 \} -\{x_1\!\to\! v_1,x_2\!\to\! v_1\}$. By Remark \ref{rem:Acyclicity}, $\cg'$ is a DAG. The new collider at $v_1$ is not a support collider by~\ref{item:T2}. So $\cg'$ has exactly two more edges of $T$ than $\cg$ and two fewer support colliders (at $u_2$ and $v_1$).

\emph{3b:} Next, suppose that $|\pa(v_1) |=2$ and $\pa(v_1) \subseteq N(v_1)$. By Lemma~\ref{lem: paths in T bar}, every directed path $u_2\rightsquigarrow v_1$ in $\cg$ has all vertices before $v_1$ in $N(v_1)$. Hence, 
\[
\text{$(\star)$ a directed path starting $u_2\to x$ with $x\notin N(v_1)\cup\{v_1\}$ never meets $N(v_1)\cup\{v_1\}$.}
\]Build $\cg'$ from $\mathcal{B}$ in the following way. 
\begin{enumerate}[label=(\roman*)]
\item For each edge $u_2\to u_i$ with $u_i\in N(v_1)$, delete it; if $u_i$ also had a parent $u_j\in N(v_1)$, add $u_a\to u_i$ for some $a\in\{1,3\}$ with $u_a\ne u_j$ so that $u_i$ remains a support collider. Since $\mathrm{indeg}(u_i) \leq 2$, at least one of $u_1$ and $u_3$ is not a parent of $u_i$. Note that $u_i\notin\{u_1,u_3\}$, since this would create a directed cycle in $\cg$.

\item If $u_2\to v_1\in E(\cg)$, delete it and add $u_b\to v_1$ for some $b\in\{1,3\}$ with $u_b\notin \pa_\cg(v_1)$. Such a $b$ exists for the same reason as in (i). 
\item Add $v_1\to u_2$.
\end{enumerate}
By construction, all of the in-degrees are at most 2. The edges of $T$ included in the skeleton of $\cg'$ increase from those of $\cg$ by exactly one ($v_1u_2$, or $v_1u_b$ if $u_2\in \pa(v_1)$). The digraph $\cg'$ has one fewer support collider than $\cg$, as the replaced ones in (i) are not already colliders, else $u_a$ would be a child of $u_i$. Before the addition of $v_1 \rightarrow u_2$ in step (iii), the graph is acyclic. For the sake of contradiction, suppose there is a cycle in $\cg'$ through $v_1\to u_2$, so that $\cg'-(v_1\to u_2)$ has a directed path $u_2\rightsquigarrow v_1$. Since we have deleted all edges $u_2 \rightarrow u_i$ with $u_i \in N(v_1)$ and the edge $u_2\to v_1$, and all added edges start in $N(v_1)$, this path starts with an edge $u_2\to x$ of $\cg$ with $x\notin N(v_1)\cup\{v_1\}$. By $(\star)$, the path then stays in $\cg$ and never meets $N(v_1)\cup\{v_1\}$, so it cannot reach $v_1$, a contradiction. So $\cg'$ is a DAG.

\emph{3c:} Suppose $|\pa(v_1)|=2$ and we are in neither Case 3a nor Case 3b. Then $v_1$ is not a support collider, and some parent $y$ of $v_1$ is not in $N(v_1)$; let $x$ be the other parent of $v_1$. If $x \ne u_2$, set $z=u_2$. If $x=u_2$, set $z=u_1$. Put $\cg'=(\mathcal{B}+z\!\to\! v_1)-\{y\!\to\! v_1\}$. Then $\cg'$ is acyclic by Remark \ref{rem:Acyclicity} and has all in-degrees at most $2$ by construction. Also, $z \rightarrow v_1$ is a new tree edge in $\cg'$ and $y \rightarrow v_1$ was not. So $\cg'$ has exactly one more edge of $T$ than $\cg$ and one fewer support collider.

\emph{3d:} Finally, suppose $|\pa(v_1)|=1$ and let $p$ be the unique parent of $v_1$. Let $z=u_1$ if $p=u_2$ and $z=u_2$ otherwise. Define $\cg'=\mathcal{B} +z\!\to\! v_1$. This is again acyclic by Remark \ref{rem:Acyclicity} and has in-degree at most 2 by construction. The DAG $\cg'$ has one more tree edge than $\cg$ and one fewer support collider.

In every case, $\ba_T \cdot \bc_{\cg'} = \ba_T \cdot \bc_{\cg}$, and $\cg'$ has fewer support colliders than $\cg$, as needed.
\end{proof}

\subsection{Facet-defining}\label{ssec:facet}
We will now prove that the matrix of vertices of $\bdcim n2$ that achieve equality on the outer inequality of $T$ has maximal rank. This combined with Proposition \ref{prop:step} will show that the inequality is facet-defining.

\begin{proposition}\label{prop: tree facet}
    Let $T$ be a tree on $m$ vertices with no degree $2$ nodes and let $M_T$ denote the set of vertices $\bc_\cg$ of $\bdcim n2$ for which $\ba_T \cdot \bc_\cg = m-1.$ Then $\rank(M_T) = \binom n2 + \binom n3$.
\end{proposition}

\begin{proof}
Again by Lemma~\ref{lem: building facets}, it suffices to treat the case $n=m$, which we assume throughout this proof. Recall that $\Lambda$ consists of the pairs that are not edges of $T$ and the triples that are not in $\Nb$. We choose tight DAGs and take their imsets as the columns of a matrix $M$, which will be a submatrix of $M_T$. For each column we give its \emph{residual}, which we define as $\bc_\cg-\mathbf 1_T$. Two pieces of notation are needed to describe the columns.
\begin{itemize}[leftmargin=*,nosep]
\item For $r\in[m]$, let $\mathcal A_r$ be the DAG obtained from $T$ by directing every edge away from $r$. Every vertex other than $r$ then has exactly one parent, so $\bc_{\mathcal A_r}=\mathbf 1_T$. Note that any choice of root $r$ yields the same imset.
\item A \emph{$2$-path of $T$ with center $k$} is a triple $\{i,j,k\}$ with $ik,jk\in E(T)$, i.e.\ the vertex set of a path $i\!-\!k\!-\!j$ with two edges. If a DAG contains both edges directed into $k$, then $\{i,j,k\}$ is a collider in that DAG.
\end{itemize}
Let $M$ be the matrix with the following columns. Each type of column is depicted in Figure \ref{fig:columns} for the tree $T_0$ on vertex set $[6]$ with edge set $E(T_0) = \{13, 23, 34, 45, 46\}$.
\begin{enumerate}[label=\textbf{(C\arabic*)},leftmargin=*,nosep]
\item The imset of the DAG $\mathcal A_r$ for a fixed root $r$; in Figure \ref{fig:columns}, we have used $r=1$. Here $\bc_{\mathcal A_r}=\mathbf 1_T$, so the residual is $0$.
\item For each $2$-path $\{i,j,k\}$ in $T$ with center $k$, the imset of $\mathcal A_k$ with the edges $ki$ and $kj$ reversed ($i\to k\leftarrow j$). The residual is $\be_{ijk}$.
\item For each non-edge $ab\notin E(T)$, let $a=x_0,\dots,x_r=b$ be the unique $ab$-path in $T$, so that $r \geq 2$; orient this path $x_0\to x_1\leftarrow x_2\leftarrow\cdots\leftarrow x_r$ and orient the other edges of $T$ away from the vertices of the path. Finally add the edge $a\to b$; this yields a DAG since the only cycle in the skeleton consists of the $ab$-path and the edges $ab$, which is not a directed cycle in the constructed digraph. The residual is $\be_{ab}+\be_{ax_1x_2}$; we note that the second term of this residual is the residual of the graph constructed in (C2) for the 2-path $a, x_1 x_2$ with center $x_1$.
\item For each triple $\tau\in\Lambda$ that does not induce a $2$-paths in $T$, we construct a DAG in the following way depending on whether there is an edge between any pair of vertices in $\tau$.
\begin{enumerate}[label = (\alph*)]
\item Suppose that the vertices of $\tau$ have exactly one $T$-edge among them, $pk$, and third vertex $q$ not adjacent to either $p$ or $k$. Let $c\in\{p,k\}$ be the node which is not in the opposite connected component from $q$ in $T-pk$; the parent of $c$ in $\mathcal A_q$ is the other endpoint $c'$ of $pk$, so $\mathcal A_q+(q\to c)$ has collider $\tau$ at $c$; residual $\be_{qc}+\be_{pkq}$. \item Now suppose $\tau=\{k,p,q\}\subseteq[m]$ has no $T$ edge among its vertices. Let $\pi$ be the parent map of $\mathcal A_k$; that is, the map which sends each vertex of $\mathcal A_k$ to its unique parent. If $\pi(p)\ne\pi(q)$ let $\mathcal O_{kpq}$ be $\mathcal A_k$ with the parent edges of $p,q$ reversed. If $\pi(p)=\pi(q)=w$ then $w\notin N(k)$ since $\tau \not \in\Nb$, so $w$ has depth at least $2$. In this case, to construct $\mathcal O_{kpq}$ we reverse only the parent edge of $w$. Since $w$ has depth at least $2$, the reversed edges do not point into $k$, so $k$ is still a source. All in-degrees remain at most 2, and the $T$-paths from $k$ to $p,q$ contain a reversed edge, so $\mathcal O_{kpq}$ has no directed path from $k$ to $p$ or $q$. In either case, $\cg=\mathcal O_{kpq}+(p\to k)+(q\to k)$ is acyclic. 

To compute its residual, note that in $\mathcal A_k$, every vertex other than $k$ has exactly one parent, and reversing the parent edge $\pi(x)\to x$ of a vertex $x$ with $\pi(x)\neq k$ makes $x$ a second parent of $\pi(x)$. The only new collider this creates is the $2$-path $\{\pi(\pi(x)),\pi(x),x\}$ with center $\pi(x)$. If $\pi(p)\neq\pi(q)$, we reverse the parent edges of $p$ and $q$; here $\pi(p),\pi(q)\neq k$ because $p,q\notin N(k)$, and since $\pi(p)\neq\pi(q)$ no vertex receives three parents. The colliders of $\mathcal O_{kpq}$ are then exactly the two $2$-paths
$
\theta_p=\{\pi(\pi(p)),\pi(p),p\}\quad\text{and}\quad\theta_q=\{\pi(\pi(q)),\pi(q),q\}.
$

If $\pi(p)=\pi(q)=w$, we reverse only the parent edge of $w$; since $w$ has depth at least $2$, $\pi(w)\neq k$, and $\mathcal O_{kpq}$ has exactly one collider, the $2$-path $\theta_w=\{\pi(\pi(w)),\pi(w),w\}$. Let $\Theta_{kpq}$ denote this set of one or two $2$-paths. Adding $p\to k$ and $q\to k$ adds the edges $pk,qk$ and the collider $p\to k\leftarrow q$. All edges of $T$ are still present, and none of the colliders lies in $\Nb$: each $2$-path contains an edge of $T$, and $\{k,p,q\}\in\Lambda$. Hence $\cg$ is tight, with residual
\[
\bc_\cg-\mathbf 1_T=\be_{pk}+\be_{qk}+\be_{kpq}+\sum_{\theta\in\Theta_{kpq}}\be_\theta .
\]
\end{enumerate}
\item \emph{Support columns:} for each internal vertex $v$, $\sigma=\{i,j,k\}\subseteq N(v)$ and each $k\in\sigma$ we obtain the following DAG: let $C_k$ be the component of $T-vk$ containing $k$; orient $C_k$ away from $k$, orient $i\to v\leftarrow j$ and the rest of the other component away from $\{i,j,v\}$, delete $vk$, and add $i\to k\leftarrow j$. This digraph is acyclic, contains $E(T)\setminus\{vk\}$ and exactly one support collider $\sigma$ (at $k$), so it is tight, with residual $\be_{ijk}-\be_{vk}+\be_{ik}+\be_{jk}+\be_{ijv}.$ This construction results in $3 | \Nb|$ columns of $M$.
\end{enumerate}
So $M$ has $1+|\Lambda|+3|\Nb|$ columns. Figure~\ref{fig:columns} shows one of each type of column for the double star $T_0$ with edges $13,23,34,45,46$. Here, $n=m=6$, $d=35$, $|S|=7$, $|\Lambda|=28$, so there are $35$ columns.
 
\begin{figure}[ht]
\centering
\begin{minipage}[t]{0.31\linewidth}\centering
\begin{tikzpicture}[scale=0.62,transform shape]
\node[vtx] (n1) at (0,0.9) {1};
\node[vtx] (n2) at (0,-0.9) {2};
\node[vtx] (n3) at (1.1,0) {3};
\node[vtx] (n4) at (2.5,0) {4};
\node[vtx] (n5) at (3.6,0.9) {5};
\node[vtx] (n6) at (3.6,-0.9) {6};
\draw[tedge] (n1) -- (n3);
\draw[tedge] (n3) -- (n2);
\draw[tedge] (n3) -- (n4);
\draw[tedge] (n4) -- (n5);
\draw[tedge] (n4) -- (n6);
\end{tikzpicture}\\[2pt]
{\footnotesize \textbf{(C1) first column: $\mathcal A_1$}\\ $c=\mathbf 1_T$}
\end{minipage}\hfill\begin{minipage}[t]{0.31\linewidth}\centering
\begin{tikzpicture}[scale=0.62,transform shape]
\node[vtx] (n1) at (0,0.9) {1};
\node[vtx] (n2) at (0,-0.9) {2};
\node[vtx,coll] (n3) at (1.1,0) {3};
\node[vtx] (n4) at (2.5,0) {4};
\node[vtx] (n5) at (3.6,0.9) {5};
\node[vtx] (n6) at (3.6,-0.9) {6};
\draw[tedge] (n1) -- (n3);
\draw[tedge] (n2) -- (n3);
\draw[tedge] (n3) -- (n4);
\draw[tedge] (n4) -- (n5);
\draw[tedge] (n4) -- (n6);
\end{tikzpicture}\\[2pt]
{\footnotesize \textbf{(C2) $2$-path $\{1,2,3\}$}\\ $c-\mathbf 1_T=$ $\be_{123}$}
\end{minipage}\hfill\begin{minipage}[t]{0.31\linewidth}\centering
\begin{tikzpicture}[scale=0.62,transform shape]
\node[vtx] (n1) at (0,0.9) {1};
\node[vtx] (n2) at (0,-0.9) {2};
\node[vtx,coll] (n3) at (1.1,0) {3};
\node[vtx] (n4) at (2.5,0) {4};
\node[vtx] (n5) at (3.6,0.9) {5};
\node[vtx] (n6) at (3.6,-0.9) {6};
\draw[tedge] (n1) -- (n3);
\draw[xedge] (n1) -- (n5);
\draw[tedge] (n3) -- (n2);
\draw[tedge] (n4) -- (n3);
\draw[tedge] (n4) -- (n6);
\draw[tedge] (n5) -- (n4);
\end{tikzpicture}\\[2pt]
{\footnotesize \textbf{(C3) non-edge $15$}\\ $c-\mathbf 1_T=$ $\be_{15}+\be_{134}$}
\end{minipage}
 
\medskip
\begin{minipage}[t]{0.31\linewidth}\centering
\begin{tikzpicture}[scale=0.62,transform shape]
\node[vtx,coll] (n1) at (0,0.9) {1};
\node[vtx] (n2) at (0,-0.9) {2};
\node[vtx] (n3) at (1.1,0) {3};
\node[vtx] (n4) at (2.5,0) {4};
\node[vtx] (n5) at (3.6,0.9) {5};
\node[vtx] (n6) at (3.6,-0.9) {6};
\draw[tedge] (n3) -- (n1);
\draw[tedge] (n3) -- (n2);
\draw[tedge] (n4) -- (n3);
\draw[tedge] (n4) -- (n6);
\draw[xedge] (n5) -- (n1);
\draw[tedge] (n5) -- (n4);
\end{tikzpicture}\\[2pt]
{\footnotesize \textbf{(C4a) one $T$-edge: $\{1,3,5\}$}\\ $c-\mathbf 1_T=$ $\be_{15}+\be_{135}$}
\end{minipage}\hfill\begin{minipage}[t]{0.31\linewidth}\centering
\begin{tikzpicture}[scale=0.62,transform shape]
\node[vtx,coll] (n1) at (0,0.9) {1};
\node[vtx] (n2) at (0,-0.9) {2};
\node[vtx,coll] (n3) at (1.1,0) {3};
\node[vtx,coll] (n4) at (2.5,0) {4};
\node[vtx] (n5) at (3.6,0.9) {5};
\node[vtx] (n6) at (3.6,-0.9) {6};
\draw[tedge] (n1) -- (n3);
\draw[xedge] (n2) -- (n1);
\draw[tedge] (n2) -- (n3);
\draw[tedge] (n3) -- (n4);
\draw[tedge] (n4) -- (n6);
\draw[xedge] (n5) -- (n1);
\draw[tedge] (n5) -- (n4);
\end{tikzpicture}\\[2pt]
{\footnotesize \textbf{(C4b) no $T$-edge: $\{1,2,5\}$}\\ $c-\mathbf 1_T=$ $\be_{12}+\be_{15}+\be_{123}+\be_{125}+\be_{345}$}
\end{minipage}\hfill\begin{minipage}[t]{0.31\linewidth}\centering
\begin{tikzpicture}[scale=0.62,transform shape]
\node[vtx] (n1) at (0,0.9) {1};
\node[vtx] (n2) at (0,-0.9) {2};
\node[vtx,coll] (n3) at (1.1,0) {3};
\node[vtx,coll] (n4) at (2.5,0) {4};
\node[vtx] (n5) at (3.6,0.9) {5};
\node[vtx] (n6) at (3.6,-0.9) {6};
\draw[tedge] (n1) -- (n3);
\draw[xedge] (n1) -- (n4);
\draw[tedge] (n2) -- (n3);
\draw[xedge] (n2) -- (n4);
\draw[tedge] (n4) -- (n5);
\draw[tedge] (n4) -- (n6);
\end{tikzpicture}\\[2pt]
{\footnotesize \textbf{(C5) $v=3,\ \sigma=\{1,2,4\},\ k=4$}\\ $c-\mathbf 1_T=$ $\be_{14}+\be_{24}-\be_{34}+\be_{123}+\be_{124}$}
\end{minipage}
\caption{One column of each type for $T_0$, with residual $\bc_\cg-\mathbf 1_T$. The shaded purple vertices are the ones with two parents. The dashed purple edges correspond to non-edges of $T_0$.}\label{fig:columns}
\end{figure}
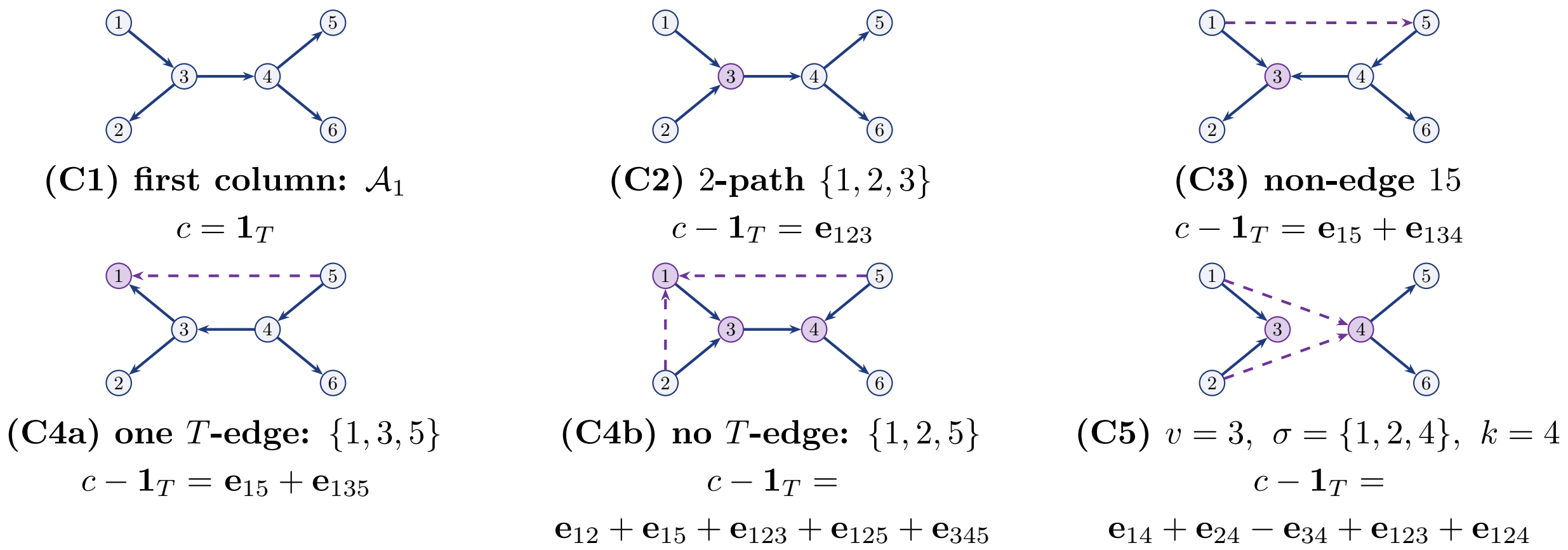
 
  We now compute the rank of $M$. Since $\bc_{\mathcal A_r}=\mathbf 1_T$ is itself a column of $M$, subtracting it from every other column does not change the rank of $M$. Let $M'$ be the matrix whose columns are the residuals $\bc_\cg-\mathbf 1_T$ of the DAGs (C2)-(C5), with rows indexed by $S\cup\Lambda$. The residual of (C1) is the zero vector, so we do not include it in $M'$, and $M'$ has $|\Lambda|+3|\Nb|$ columns. Then, $\rank(M)=\rank\bigl(\mathbf 1_T\mid M'\bigr)$, and we will show that $\rank M'= d-1$ and that $\mathbf 1_T$ is not in the column span of $M'$. The proof that $\rank M'=d-1$ consists of three parts: the enumeration of the entries of $M'$; the reduction of $M'$ to a block form; and the computation of the rank of the blocks. The proof consists of three parts: the enumeration of the entries of $M'$; the reduction of $M'$ to a block form; and the computation of the rank of the blocks.

To order the rows we split $\Lambda$ into three \emph{levels}: $\Lambda_1$, the $2$-paths; $\Lambda_2$, the pairs that are not edges of $T$; $\Lambda_3$, the triples of $\Lambda$ that do not induce $2$-paths in $T$. Here $\Lambda_1\subseteq\Lambda$, since a $2$-path is never in $\Nb$: the triples in $\Nb$ contain no $T$-edge, because $T$ has no triangle. We list the elements of $S$ first and those of $\Lambda$ afterwards, and we list $\Lambda$ by increasing level: first $\Lambda_1$, then $\Lambda_2$, then $\Lambda_3$. We call the \emph{constructing coordinate} of a column of type (C2), (C3) or (C4) the $2$-path $\{i,j,k\}$, the pair $ab$, or the triple $\tau$, respectively, from which the column is built. Table~\ref{tab:land} lists the residuals computed above, sorted by the level of each coordinate; this sorting is immediate from the definitions of $\Lambda_1$, $\Lambda_2$ and $\Lambda_3$.
 
\begin{table}[ht]\centering\small
\begin{tabular}{@{}lllp{6.4cm}@{}}
\toprule
column & number & $S$-part of $c-\mathbf 1_T$ & $\Lambda$-part of $c-\mathbf 1_T$\\
\midrule
(C1) $\mathcal A_r$ & $1$ & $0$ & $0$ (not a column of $M'$)\\
(C2) $2$-path $ijk$ & $|\Lambda_1|$ & $0$ & $\be_{ijk}$ (constructing coordinate)\\
(C3) non-edge $ab$ & $|\Lambda_2|$ & $0$ & $\be_{ab}$ (constr.) $+\;\be_{x_0x_1x_2}$ (level 1)\\
(C4a) triple, one $T$-edge & \multirow{2}{*}{$|\Lambda_3|$} & $0$ & $\be_{pkq}$ (constr.) $+\;\be_{qc}$ (level 2)\\
(C4b) triple, no $T$-edge & & $0$ & $\be_{kpq}$ (constr.) $+\;\be_{pk}+\be_{qk}$ (level 2) $+\sum_{\theta\in\Theta(\mathcal O)}\be_\theta$ (level 1)\\
(C5) support & $3|\Nb|$ & $\be_{ijk}-\be_{vk}$ & $\be_{ik}+\be_{jk}+\be_{ijv}$ (no constr. coordinate)\\
\bottomrule
\end{tabular}
\caption{The entries of the residuals of the columns. We use the abbreviation ``constr." to indicate that a coordinate in the support of the residual is the constructing coordinate of the column.}\label{tab:land}
\end{table}
 
The columns (C2)-(C4) are in bijection with $\Lambda$ via their constructing coordinates: every $2$-path occurs in (C2), every non-edge in (C3), and every triple of $\Lambda$ that is not a $2$-path in (C4). By Table~\ref{tab:land} such a column has zero $S$-part, entry $1$ in its constructing coordinate, and otherwise only entries equal to $1$ in coordinates of strictly lower level (the $2$-paths $\theta\in\Theta(\mathcal O)$ in (C4)(b) are distinct, since a $2$-path has a unique center). Consequently the $\Lambda\times\Lambda$ block $B$ formed by these columns is upper triangular with ones on the diagonal, and
\[\rank(B)=|\Lambda|.\]

The columns (C5) have the $S$-part $\be_{\sigma}-\be_{vk}$ and a $\Lambda$-part, which we collect in a matrix $W$. Ordering the columns as (C5) first and then (C2)--(C4), and the rows as $S$ first and then $\Lambda$, we obtain
\[
M'=\begin{pmatrix}U&0\\ W&B\end{pmatrix},
\]
where $U$ is the $S\times 3|\Nb|$ matrix with columns $\be_{ijk}-\be_{vk}$. Since the upper right block is zero,
\[
\rank(M')\geq\rank(U)+\rank(B)=\rank(U)+|\Lambda|.\tag{3}
\]
Indeed, if a linear combination of linearly independent columns of $U$ and of $B$, viewed as columns of $M'$, vanishes, then its restriction to the rows indexed by $S$ involves only the columns of $U$, so their coefficients vanish, and then the coefficients of the columns of $B$ vanish as well.
 
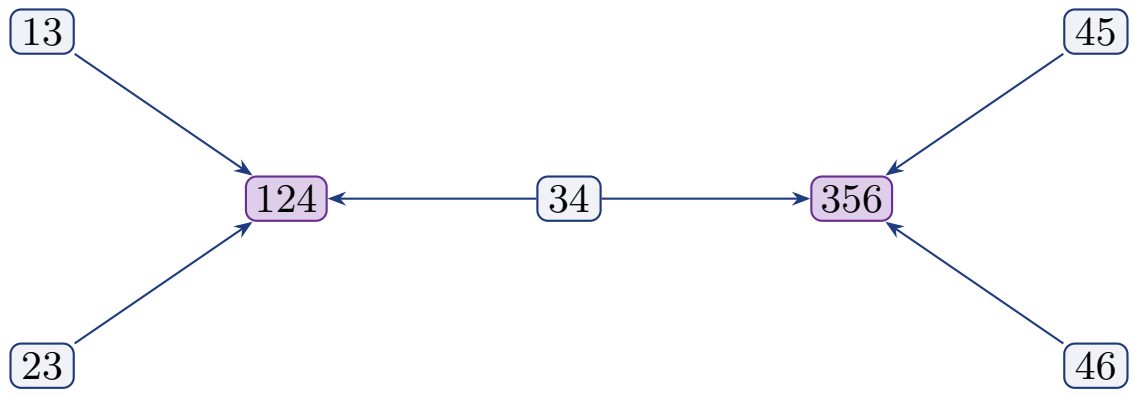
\begin{figure}[ht]\centering
\begin{tikzpicture}[scale=0.85,transform shape]
\node[hvtx] (h13) at (0,1.3) {\scriptsize $13$};
\node[hvtx] (h23) at (0,-1.3) {\scriptsize $23$};
\node[hvtx,coll] (h124) at (1.9,0) {\scriptsize $124$};
\node[hvtx] (h34) at (4.1,0) {\scriptsize $34$};
\node[hvtx,coll] (h356) at (6.3,0) {\scriptsize $356$};
\node[hvtx] (h45) at (8.2,1.3) {\scriptsize $45$};
\node[hvtx] (h46) at (8.2,-1.3) {\scriptsize $46$};
\draw[harc] (h13) -- (h124);
\draw[harc] (h23) -- (h124);
\draw[harc] (h34) -- (h124);
\draw[harc] (h34) -- (h356);
\draw[harc] (h45) -- (h356);
\draw[harc] (h46) -- (h356);
\end{tikzpicture}
\caption{The DAG $\ch$ for $T_0$ (purple: the triples $124,356$). The two stars of $T_0$ meet in the vertex $34$, so $\ch$ is connected and $\rank U=6$; hence $\dim F=28+6=34=d-1$.}\label{fig:H}
\end{figure}

It remains to compute $\rank (U)$. Each column of $U$ is $\be_{\sigma}-\be_{vk}$, that is, the incidence vector of the edge $vk\to\sigma$. Hence $U$ is the incidence matrix of the DAG $\ch$ on the vertex set $S=E(T)\cup\Nb$ with an edge $vk\to\sigma$ for every internal vertex $v$, every $\sigma\subseteq N(v)$ with $|\sigma|=3$, and every $k\in\sigma$ (see Figure~\ref{fig:H} for $T_0$). We show that $\ch$ is connected when the directions of the edges are ignored.
 
Let $\sigma\in\Nb$. There is a unique vertex $v$ with $\sigma\subseteq N(v)$, since two common neighbors of three vertices would produce a cycle in $T$, and the edges of $\ch$ entering $\sigma$ are exactly the three edges $vi\to\sigma$, $vj\to\sigma$, $vk\to\sigma$ for $\sigma=\{i,j,k\}$. In particular every $\sigma\in\Nb$ is joined to an edge of $T$ in $\ch$, so it suffices to show that any two edges of $T$ are joined by a walk in $\ch$.
 
First, let $v$ be a vertex of $T$ that is not a leaf and let $vk$, $vk'$ be two distinct edges at $v$. By the assumption of the theorem that $T$ has no vertex of degree $2$, we have $\deg v\ge3$; this is the point at which this assumption enters the argument. Hence, there is a vertex $k''\in N(v)\setminus\{k,k'\}$, and $\sigma=\{k,k',k''\}\in\Nb$. The edges $vk\to\sigma$ and $vk'\to\sigma$ show that $vk,\sigma,vk'$ is a path in $\ch$.
 
Second, let $e,f\in E(T)$ and let $P$ be the shortest path in $T$ that contains $e$ and $f$. Denote its edges, in order, by $e=e_0,e_1,\dots,e_s=f$, and let $w_t$ be the common vertex of $e_{t-1}$ and $e_t$ for $t=1,\dots,s$. The vertex $w_t$ lies on two edges, so it is not a leaf. By the first step, which uses once more that $T$ has no vertex of degree $2$, there is a triple $\sigma_t\in\Nb$ such that $e_{t-1},\sigma_t,e_t$ is a walk in $\ch$. Concatenating these paths gives the walk $e_0,\sigma_1,e_1,\sigma_2,\dots,\sigma_s,e_s$ from $e$ to $f$ in $\ch$. For example, consider $T_0$ and the edges $13$ and $45$; this is the walk $13,124,34,356,45$. The graph $\ch$ for $T_0$ is pictured in Figure \ref{fig:H}.
Note that if $T$ has a single edge, then $\Nb=\emptyset$ and $\ch$ consists of one vertex.

Thus $\ch$ is connected.  By Theorem~\ref{thm:incidence matrix}, the incidence matrix of a connected directed graph with $N$ vertices has rank $N-1$, so $\rank U=|S|-1$. By (3),
\[
\rank (M')\geq|\Lambda|+|S|-1=d-1.
\]
Every column of $M$ is the imset of a tight DAG, so for every column $\bc_\cg$ of $M'$'s construction we have $\ba_T\cdot(\bc_\cg-\mathbf 1_T)=(m-1)-(m-1)=0$. Hence the columns of $M'$ lie in the hyperplane $\{\bx:\ba_T\cdot\bx=0\}$; in particular $\rank(M')=d-1$. On the other hand $\ba_T\cdot\mathbf 1_T=m-1\neq0$, so $\mathbf 1_T$ is not in the column span of $M'$. Therefore
\[
\rank(M_T)\;\geq\;\rank(M)\;=\;\rank\bigl(\mathbf 1_T\mid M'\bigr)\;=\;\rank(M')+1\;=\;d,
\]
and since $M_T$ has $d$ rows, $\rank(M_T)=\binom m2+\binom m3$.
\end{proof}

\begin{remark}
The assumption that $T$ has no vertex of degree $2$ was used in the connectedness of $\ch$ above. If $T$ has a vertex $w$ of degree $2$, no triple of $\Nb$ is centered at $w$, the edges of $T$ on the two sides of $w$ are not joined by any walk in $\ch$, and $\ch$ is disconnected. The inequality still defines a face of $\bdcim n2$ in this case, but this face has dimension smaller than $d-1$, so it is not a facet.
\end{remark}
 
\begin{proof}[Proof of Theorem \ref{thm:tree-facet}]
    Repeatedly applying Proposition \ref{prop:step} to a given DAG, we reach a DAG with no support colliders whose value \eqref{eq:value} is the number of edges of $T$ it contains, at most $m-1$. So Proposition \ref{prop:step} implies that the outer inequality associated to $T$ is valid on $\bdcim{m}{2}$. Together with Proposition \ref{prop: tree facet}, this shows that the outer inequality associated to $T$ defines a facet of $\bdcim{m}{2}$. Furthermore, by Lemma \ref{lem: building facets}, we conclude that the outer inequality defines a facet of $\bdcim{n}{2}$, for any $n\geq m$.
\end{proof}

\section{Inequalities from Cliques and Disconnected Graphs}
\label{sec:clique-ineqs}
\noindent In the previous sections, we focused on constructing facet-defining inequalities of $\bdcim{n}{2}$ whose edge variable support was given by a connected graph. In this section, we introduce several new families of facets that are supported on the edges of an undirected graph which is not necessarily connected. We also introduce more facets whose edge variables correspond to a clique in addition to the ones identified in Section \ref{sec:path-ineqs}. We first consider an inequality that arises from a $4$-clique, which is different from the one that we studied in Example \ref{eg: chordal example}.  

\begin{example}\label{eg: four clique other example}
    We consider the complete graph on vertices $1,2,3,$ and $4$. Direct computation shows that the following inequality defines a facet of $\bdcim{4}{2}$:
    $$x_{12}+x_{13}+x_{14}+x_{23}+x_{24}+x_{34} - (x_{123}+x_{124}+x_{134}+x_{234})\leq 3.$$

  For now, we assume that the inequality defines a face $F$ of $\bdcim{4}{2}$. To show that $F$ is a facet, we find several vertices of $F$, put them as columns in a matrix $M$ and compute the rank of $M$.

 The vertices we consider are the indicator vectors of the spanning trees of $K_4$. For a spanning tree $T$ of $K_4$, let $\be_T=\sum_{ij\in E(T)}\be_{ij}$. This is the characteristic imset of $T$ with its edges directed away from any root, since every vertex other than the root then has exactly one parent. Three edges of $K_4$ constitute a spanning tree exactly when they do not form a triangle. There are $20-4=16$ such columns and clearly they form a submatrix of rank $6$. In fact, we can perform simple elementary column operations to transform these columns into the identity matrix of rank $6$. 
    Then, for each triple $ijk$, with $l$ the fourth vertex, we take the DAG with edges $j\to i$, $k\to i$, $j\to k$ and $k\to l$. It has four edges and exactly one collider, at $i$, so its imset $\be_{ij}+\be_{ik}+\be_{jk}+\be_{kl}+\be_{ijk}$ lies on $F$. Subtracting suitable combinations of the tree columns, which span all edge coordinates, turns these four columns into $\be_{123},\be_{124},\be_{134},\be_{234}$.
  
    Thus, $M$ has rank $10=\dim(\bdcim{4}{2})$. Since $F$ does not contain the origin, the affine hull of $F$ has dimension $9$. This implies that $F$ is a facet of $\bdcim{4}{2}$, and by Lemma \ref{lem: building facets}, also a facet of $\bdcim{n}{2}$, for each $n\geq4$. 
\end{example}

 Once again, the previous inequality is part of a larger family that holds for $\bdcim{n}{2}$. The following counting lemma will be a key tool for proving that such facets are valid.

\begin{lemma}
\label{lemma:count-coord-sums}
Let $\cg = ([n], E)\in\dags{n}{2}$. Let $V_\cg(d) = |\{v \in [n] ~:~ |\pa_\mathcal{G}(v)| = d\}|$ be the number of vertices of $\cg$ with exactly $d$ parents. Then
\begin{align*}
\sum_{i < j} \bc_\cg(ij) = |E| = V_\cg(1) + 2 V_\cg(2),~~ 
\sum_{i < j < k} \bc_\cg(ijk) = V_\cg(2).
\end{align*}
\end{lemma}

\begin{proof}
    For the first equality, observe that
    \[
        \sum_{i < j} \bc_\cg(ij) = \sum_{i \in V(\cg)} | \pa_\cg(i) |.
    \]
    A vertex $i$ contributes $0$ to this sum if it is a source, $1$ if it has exactly one parent, and $2$ if it has exactly two parents.

    For the second equality, recall that since $|\pa_\cg(j)| \leq 2$ for all $j$, each vertex is a collider in at most one triple of vertices. So this sum counts the number of vertices that are colliders, that is, the number of vertices with exactly two parents.
\end{proof}

The following facet-defining inequality is a generalization of the one presented in Example \ref{eg: four clique other example}. It is another example of a facet whose edge-support corresponds to the edges of a clique on $m$ nodes; in this case, all edge variables have coefficient $1$. 

\begin{figure}[ht]
\centering
\setlength{\tabcolsep}{5pt}
\begin{tabular}{ccc}
\begin{tikzpicture}[>=stealth,x=0.85cm,y=0.85cm,every node/.style={circle,draw=cimnode,fill=cimnodefill,inner sep=0.6pt,minimum size=4.4mm,font=\scriptsize}]
\node (n1) at (0,1.2) {1};
\node (n2) at (1.2,1.2) {2};
\node (n3) at (0,0) {3};
\node (n4) at (1.2,0) {4};
\draw[->,thick,cimedge] (n1) -- (n2);
\draw[->,thick,cimedge] (n1) -- (n3);
\draw[->,thick,cimedge] (n1) -- (n4);
\end{tikzpicture} &
\begin{tikzpicture}[>=stealth,x=0.85cm,y=0.85cm,every node/.style={circle,draw=cimnode,fill=cimnodefill,inner sep=0.6pt,minimum size=4.4mm,font=\scriptsize}]
\node (n1) at (0,1.2) {1};
\node (n2) at (1.2,1.2) {2};
\node (n3) at (0,0) {3};
\node (n4) at (1.2,0) {4};
\draw[->,thick,cimedge] (n1) -- (n2);
\draw[->,thick,cimedge] (n2) -- (n3);
\draw[->,thick,cimedge] (n3) -- (n4);
\end{tikzpicture} &
\begin{tikzpicture}[>=stealth,x=0.85cm,y=0.85cm,every node/.style={circle,draw=cimnode,fill=cimnodefill,inner sep=0.6pt,minimum size=4.4mm,font=\scriptsize}]
\node (n1) at (0,1.2) {1};
\node (n2) at (1.2,1.2) {2};
\node[coll] (n3) at (0,0) {3};
\node (n4) at (1.2,0) {4};
\draw[->,thick,cimedge] (n1) -- (n2);
\draw[->,thick,cimedge] (n1) -- (n3);
\draw[->,thick,cimedge] (n2) -- (n3);
\draw[->,thick,cimedge] (n1) -- (n4);
\end{tikzpicture} \\
{\scriptsize star} & {\scriptsize path} & {\scriptsize triangle + pendant} \\[4pt]
\end{tabular}\\
\caption{Columns of $M$ for $m=4$: oriented spanning trees (first two) and the DAG $D_{123}$ (last), whose imset is $\be_{12}+\be_{13}+\be_{14}+\be_{23}+\be_{123}$.}
\label{fig:clique-cols}
\end{figure}
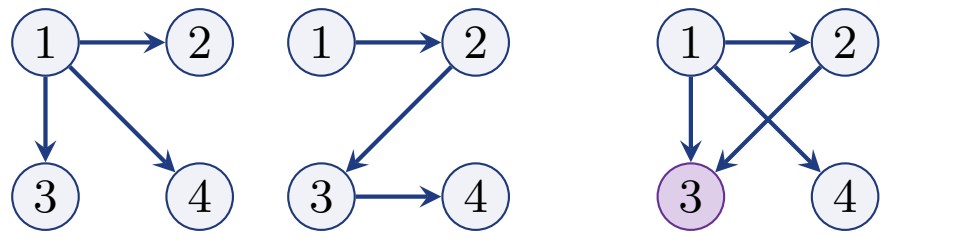

\begin{theorem}\label{thm: other clique inequality}
    The following inequality defines a facet of $\bdcim{n}{2}$ for any $n\geq m$:
    $$\sum_{1\leq i<j\leq m} x_{ij} - \sum_{1\leq i<j<k\leq m}x_{ijk} \leq m-1.$$
\end{theorem}
\begin{proof}
    As usual, we show that the above inequality defines a facet of $\bdcim{m}{2}$; hence, by Lemma \ref{lem: building facets}, it also defines a facet of $\bdcim{n}{2}$, for any $n\geq m$. 

   \textbf{Validity:} Let $\cg\in\dags{m}{2}$. Then, by Lemma \ref{lemma:count-coord-sums}
    \[
    \sum_{1\leq i<j\leq m} \bc_\cg(ij) - \sum_{1\leq i<j<k\leq m}\bc_\cg(ijk) = V_\cg(1) + 2V_\cg(2) - V_\cg(2) = V_\cg(1) + V_\cg(2). 
    \]
    The quantity $V_\cg(1) + V_\cg(2)$ is simply the number of nodes of $\cg$ which are not sources. This is at most $m-1$ since a DAG on $m$ nodes has at least one source.

  \textbf{Facet-defining:} Let $F$ be the face of $\bdcim{m}{2}$ defined by this inequality. By the above argument, a DAG $\cg$ on $[m]$ attains $\sum x_{ij}-\sum x_{ijk}=m-1$ exactly when $\cg$ has a unique source. Since the zero-vector does not belong to $F$, it suffices to exhibit such DAGs whose imsets span $\mathbb R^{\binom m2+\binom m3}$. Let $M$ be the matrix whose columns are the following imsets; see Figure~\ref{fig:clique-cols} for the $m=4$ case.
    \begin{enumerate}
        \item For every spanning tree $T$ of $K_m$, the vector $\be_T=\sum_{ij\in E(T)}\be_{ij}$. This is the characteristic imset of $T$ with its edges directed away from any chosen root.
        \item For every $i<j<k$, the imset of the DAG $D_{ijk}$ with edges $i\to j$, $i\to k$, $j\to k$ and $i\to\ell$ for all $\ell\notin\{i,j,k\}$; that is, $\sum_{\ell\neq i}\be_{i\ell}+\be_{jk}+\be_{ijk}$.
    \end{enumerate}
    First we claim that the columns of type (1) span all edge coordinates; that is, after deleting the $\binom{m}{3}$ rows indexed by collider variables, on which they vanish, they have full rank. Let $e = ab$ and $f=xy$ be distinct edges of $K_m$ and let $P$ be a path in $K_m$ from $x$ to $y$ through the edge $ab$. Extend $P$ to a spanning tree $T$ of $K_m$. Then $T-e+f$ is also a spanning tree of $K_m$. 
    Then $\be_{T} - \be_{T - e +f} = \be_e - \be_f$. These span the $(\binom{m}{2}-1)$-dimensional subspace of $\mathbb{R}^{\binom{m}{2}}$ where the sum of all coordinates is $0$. Since $\be_T$ has coordinate sum $m-1 \neq 0$ for every spanning tree $T$, the columns of type (1) span all edge coordinates and $\be_{ij}$ belongs to their span for all $ij$.  Finally, modulo the edge coordinates, the imset of $D_{ijk}$ is the unit vector $\be_{ijk}$. So, $M$ has rank $\binom m2+\binom m3=\dim\bdcim{m}{2}$ and thus, $F$ is a facet.
\end{proof}
    
We now introduce two further families of facet-defining inequalities. To prove that they are actually valid on $\bdcim{n}{2}$, we will first need the following lemma which generalizes Lemma \ref{lemma:count-coord-sums}. 

\begin{lemma}
\label{lemma:weighted-coord-count}
Let $\cg\in\dags{n}{2}$ and $\ell(\bx) = \sum_{i < j} a_{ij} x_{ij} + \sum_{i < j < k} b_{ijk} x_{ijk}$ be a linear functional. Then $\ell(\bc_\cg) = \sum_{v \in [n]} \ell_v( \pa_\cg(v))$ where for any $S \subset [n]$ of cardinality at most $2$,
\[
\ell_v(S) = \sum_{u \in S} a_{uv} + 
\begin{cases}
b_{\{v\} \cup S}, ~~  |S| = 2,\\
0,~~ |S| \leq 1.
\end{cases}
\]
\end{lemma}
\begin{proof}
    By linearity, it suffices to account for each coordinate of $\bc_\cg$ that equals $1$. A pair $\{u,v\}$ has $\bc_\cg(uv)=1$ exactly when it is an edge of $\cg$. We assign this edge to its head. If $u\to v$, then $u\in\pa_\cg(v)$, and the term $a_{uv}$ appears in $\ell_v(\pa_\cg(v))$ and in no other $\ell_w$. A triple has $\bc_\cg(ijk)=1$ exactly when some vertex $v$ of the triple has the other two as parents. Since $|\pa_\cg(v)|\le2$, this means that the triple is $\{v\}\cup\pa_\cg(v)$ with $|\pa_\cg(v)|=2$. This vertex $v$ is unique, as if two vertices of the triple each had the other two as parents, they would be parents of each other, which is impossible in a DAG. So the term $b_{\{v\}\cup\pa_\cg(v)}$ appears in $\ell_v(\pa_\cg(v))$ and in no other $\ell_w$. Summing over all $v\in[n]$ gives $\ell(\bc_\cg)$.
\end{proof}

This lemma captures the fact that the support of $\bc_\cg$ in the coordinates $x_{ij}$ is limited to those vertices which are involved in an edge, while the support on the coordinates $x_{ijk}$ simply corresponds to the nodes that have two parents. In the following results, we call $\ell_v(\pa_\cg(v))$ the \emph{contribution of $v$ to $\ell(\bc_\cg)$}. Recall that for a graph $G$, $T(G)$ refers to the set of triples of vertices which induce triangles in $G$.
    
\begin{proposition}\label{prop: path plus vertex}
    Let $G$ be a connected undirected graph on at least $3$ vertices such that $T(G)\neq\emptyset$, and let $p$ be a vertex with $p\notin V(G)$. Then the following inequality defines a facet of $\bdcim{n}{2}$ for every $n\geq |V(G)|+1$:
    \begin{equation}\label{eq: path plus vertex}
    -\sum_{(i,j)\in E(G)}x_{ij} + 2\sum_{(i,j,k)\in T(G)}x_{ijk} + \sum_{(i,j,k)\in C(G)\setminus T(G)}x_{ijk} + \sum_{(i,j)\in E(G)}x_{ijp} \leq 1.\end{equation}
\end{proposition}
\begin{proof}
    Let $m=|V(G)|+1$. As before, it suffices by Lemma~\ref{lem: building facets} to show that the inequality defines a facet of $\bdcim{m}{2}$, where we take the vertex set to be $V(G)\cup\{p\}$. Write the left-hand side of \eqref{eq: path plus vertex} as $\ell(\bx)$, where
    \[
        \ell(\bx) \;=\; \ell_G(\bx) + \sum_{ij\in E(G)}x_{ijp},
        \qquad
        \ell_G(\bx) \;=\; -\sum_{ij\in E(G)}x_{ij} + 2\sum_{ijk\in T(G)}x_{ijk} + \sum_{ijk\in C(G)\setminus T(G)}x_{ijk},
    \]
    so that $\ell_G(\bx)\le 0$ is the inequality of $G$ that defines a facet as specified in Theorem~\ref{thm:conn graph}. For a coordinate $S$ we write $\ell_S$ for the coefficient of $x_S$ in $\ell$.

    \textbf{Validity:} Let $\ch\in\dags{m}{2}$, with vertex set $V(G)\cup\{p\}$. By Lemma~\ref{lemma:weighted-coord-count}, 
    \[\ell(\bc_\ch)=\sum_{w \in V(G) \cup \{p\}}\ell_w(\pa_\ch(w)),\]
    where $\ell_w (\emptyset) = 0$, $\ell_w(\{s\}) = \ell_{sw}$, 
    and $\ell_w(\{s,t\})=\ell_{sw}+\ell_{tw}+\ell_{stw}$ if $s$ and $t$ are distinct. We bound each contribution.

    Since $\ell_{ip}=0$ for all $i$,  $\ell_{ijp}=1$ if $ij\in E(G)$, and $\ell_{ijp}=0$ otherwise, we have
    \[
        \ell_p(\pa_\ch(p)) =
        \begin{cases}
        1, & \pa_\ch(p)=\{i,j\}\text{ with }ij\in E(G), \\
        0, & \text{otherwise}.
        \end{cases}
    \]

    Now let $v\in V(G)$ and let $S\subseteq (V(G)\cup\{p\})\setminus\{v\}$ with $|S|\le 2$. Note that $\ell_v(S)$ depends only on $S$, and we will use the following bound both for $S=\pa_\ch(v)$ here and for the parent sets of the DAGs in the facet part below. If $S=\{p\}$, then $\ell_v(S)=\ell_{pv}=0$. If $S=\{u,p\}$, then $\ell_{uvp}=-\ell_{uv}$, so $\ell_v(S)=\ell_{uv}+\ell_{pv}+\ell_{uvp}=0$. If $p\notin S$, then $\ell_v(S)$ is the contribution of $v$ to the path inequality associated to $G$ as defined in Theorem \ref{thm:conn graph}, which we claim is at most $0$. Indeed, this is trivial if $|S|\le 1$. If $S=\{u,w\}$, then $\ell_v(S)=\ell_{uv}+\ell_{wv}+\ell_{uvw}$, which is $-1-1+2=0$ when $uvw\in T(G)$, is $-1$ or $0$ when $uvw\in C(G)\setminus T(G)$, and is $\ell_{uv}+\ell_{wv}\le 0$ otherwise. So $\ell_v(S)\le 0$ for every $v\in V(G)$ and every such $S$; applied to $S=\pa_\ch(v)$, every vertex of $G$ contributes at most $0$ and $p$ contributes at most $1$, which gives $\ell(\bc_\ch)\le 1$.

    \textbf{Facet-defining:} Let $F$ be the face of $\bdcim{m}{2}$ defined by \eqref{eq: path plus vertex}. Since $F$ does not contain the origin, it suffices to find vertices of $F$ whose imsets span $\mathbb R^{\binom m2+\binom m3}$; their affine hull is then contained in the hyperplane $\ell(\bx)=1$ and has dimension $\binom m2+\binom m3-1$.  We have $\ell(\bc_\ch)=1$ if and only if the two parents of $p$ form an edge of $G$ and every $v\in V(G)$ contributes $0$. In particular, if $\ch'$ is a DAG on $V(G)$ that attains equality in the path inequality, $\ell_G(\bc_{\ch'})=0$, then adding the vertex $p$ with two parents that are adjacent in $G$ gives a DAG whose imset lies on $F$. We split the coordinates of $\mathbb{R}^{\binom{m}{2} + \binom{m}{3}}$ into the \emph{$G$-coordinates}, indexed by pairs and triples of vertices in $V(G)$, and the \emph{$p$-coordinates}, indexed by pairs $ip$ and triples $ijp$ with $i,j\in V(G)$. For an edge $kl\in E(G)$ let $\ch_{kl}$ be the DAG with edges $k\to p\leftarrow l$, with all other vertices isolated. We use the imsets of the following DAGs. In each of them $p$ has the parents $k,l$ with $kl\in E(G)$. It is straightforward to check using an analogous argument to the proof of validity that for each of these DAGs, every vertex of $G$ contributes $0$ to $\ell$; so all of their imsets lie on $F$.
    \begin{enumerate}[label=\textbf{(C\arabic*)}]
        \item $\ch_{kl}$ for every $kl\in E(G)$. Its imset is $\bc_{\ch_{kl}} = \be_{kp}+\be_{lp}+\be_{klp}$.
        \item For every $j\in V(G)$, the DAG $\ch_{kl}+(p\to j)$ for one edge $kl\in E(G)$ with $j\notin\{k,l\}$. Such an edge exists because $G$ contains a triangle: if $j$ lies on a triangle, take the edge of the triangle opposite $j$, and otherwise take any edge of the triangle. The imset is  $\bc_{\ch_{kl}}+\be_{jp}$.
        \item For each pair of distinct $i,j\in V(G)$, we include the DAG $\ch_{kl}+(p\to j)+(i\to j)$, where $kl$ is again an edge of $G$ with $j \neq k,l$ as in (C2). Here $j$ has the parents $p$ and $i$, so the imset is $\bc_{\ch_{kl}}+\be_{jp} + \be_{ij}+\be_{ijp}$.
        \item Fix a single edge $k_0 l_0 \in E(G)$.
        Let $\ch'$ be a DAG on $V(G)$ in which every vertex has at most two parents, such that $\bc_{\ch'}$ lies on the facet of $\bdcim{m-1}{2}$ defined by the path inequality of $G$ as described in Theorem~\ref{thm:conn graph}. For each such DAG, we include the DAG on $m$ vertices, $\ch' + \ch_{k_0 l_0}$, which has imset $\bc_{\ch'}+\be_{k_0p}+\be_{l_0p}+\be_{k_0l_0p}$.
    \end{enumerate}

    Let $W$ be the span of the imsets of these DAGs. We will show that $W=\mathbb R^{\binom m2+\binom m3}$. First, for every vertex $\bc_{\ch'}$ of the path facet, subtracting the column (C1) for the edge $k_0l_0$, that is $\bc_{\ch_{k_0 l_0}}$, from the corresponding column (C4) gives $\bc_{\ch'}$, and so $\bc_{\ch'}\in W$. The path inequality of $G$ defines a facet of $\bdcim{m-1}{2}$, the polytope of DAGs on the $m-1$ vertices of $V(G)$ (Theorem~\ref{thm:conn graph}). This facet contains the origin and has dimension $\binom{m-1}2+\binom{m-1}3-1$, so its vertices span the hyperplane $\{\ell_G=0\}$ of $\mathbb R^{\binom{m-1}2+\binom{m-1}3}$. We regard each $\bc_{\ch'}$ as a vector of $\mathbb R^{\binom m2+\binom m3}$ whose $p$-coordinates are all $0$. Hence $W$ contains the linear subspace
    \[
    H_G=\Bigl\{\bx\in\mathbb R^{\binom m2+\binom m3}:\ \ell_G(\bx)=0\ \text{ and }\ x_S=0\ \text{for every $p$-coordinate }S\Bigr\},
    \]
    which has dimension $\binom{m-1}2+\binom{m-1}3-1$.
    
    Next, let $j\in V(G)$ and choose an edge $kl\in E(G)$ with $j\notin\{k,l\}$. Subtracting the imset of $\ch_{kl}$ in (C1) from the imset of $\ch_{kl}+(p\rightarrow j)$ in (C2) gives $\be_{jp}$. Thus, $\be_{jp}\in W$ for every $j\in V(G)$. In particular $\be_{kp},\be_{lp}\in W$ for every edge $kl\in E(G)$, and then the column (C1) for $kl$ gives $\be_{klp}\in W$. Similarly, for distinct $i,j\in V(G)$ and an edge $kl\in E(G)$ with $j\notin\{k,l\}$, subtracting the imset of $\ch_{kl}+(p\rightarrow j)$ in (C2) from the imset of $\ch_{kl}+(p\rightarrow j)+(i\rightarrow j)$ in (C3) gives $\be_{ij}+\be_{ijp}\in W$. Now let $i,j\in V(G)$ be distinct. If $ij\notin E(G)$, then $\ell_G(\be_{ij})=0$, so $\be_{ij}\in W$, and therefore $\be_{ijp}=(\be_{ij}+\be_{ijp})-\be_{ij}\in W$. If $ij\in E(G)$, then $\be_{ijp}\in W$ as shown above, and therefore $\be_{ij}=(\be_{ij}+\be_{ijp})-\be_{ijp}\in W$. Hence, $W$ contains the unit vector of every coordinate involving $p$. 
    
    Finally, $W$ contains $H_G$ together with the vector $\be_{ij}$ for an edge $ij\in E(G)$. This vector has all $p$-coordinates equal to $0$ and satisfies $\ell_G(\be_{ij})=-1\neq 0$, so $H_G+\mathbb R\be_{ij}$ is the whole space of vectors whose $p$-coordinates vanish. Hence $W$ contains the unit vector of every coordinate indexed by a subset of $V(G)$. Together with the unit vectors of the $p$-coordinates, this gives $W=\mathbb R^{\binom m2+\binom m3}$, and $F$ is a facet.
\end{proof}

The next family combines the clique inequality of Theorem \ref{thm: other clique inequality} with isolated vertices. Its edge variables are exactly the edges of the undirected graph $K_q\sqcup \overline{K_p}$, the disjoint union of a clique on $v_1,\dots,v_q$ and $p$ isolated vertices $u_1,\dots,u_p$. The isolated vertices appear only through collider variables. As for the inequality of Proposition~\ref{prop: path plus vertex}, whose edge support is the connected graph $G$ together with the isolated vertex $p$, this is a facet whose edge support is a disconnected graph, unlike the inequalities in the previous sections.

\begin{theorem}
    Let $p,q\geq 2$, and let $v_1, \dots, v_q,u_1,\dots,u_p$ be distinct vertices. Then, the following inequality defines a facet of $\bdcim{n}{2}$ for $n\geq p+q$:
    $$\sum_{i<j\in [q]}x_{v_i v_j} - \sum_{i<j<k\in [q]}x_{v_i v_j v_k} + \sum_{i\in [q], j<k\in[p]}x_{v_iu_ju_k} + \sum_{i < j < k \in [p]} x_{u_i u_j u_k}  \leq p+q-2.$$ 
\end{theorem}
\begin{proof}
        \textbf{Validity:} Again, we let $\ell$ denote the linear functional on the left side of the inequality, and write $V=\{v_1,\dots,v_q\}$, $U=\{u_1,\dots,u_p\}$ and $W=V\cup U$. Every variable of $\ell$ is indexed by a subset of $W$, and the restriction of the imset of a DAG to the coordinates indexed by subsets of $W$ is the imset of its induced subgraph on $W$. So it suffices to prove validity for DAGs on the vertex set $W$, and by Lemma~\ref{lem: building facets} it also suffices to prove that the inequality is facet-defining for $n=p+q$. Hence, from now on we consider DAGs on the vertex set $W$. Let $\cg\in\dags{p+q}{2}$ be such a DAG. By Lemma~\ref{lemma:weighted-coord-count}, $\ell(\bc_\cg)=\sum_{w\in W}\ell_w(\pa_\cg(w))$. For $w=v_i$, every parent in $V$ contributes $1$ through an edge variable, and a pair of parents contributes $-1$ if both are in $V$, $+1$ if both are in $U$ and $0$ otherwise. By examining all cases, we see that
    \begin{equation}\label{eq: clique and independent set contributions}
    \ell_{v_i}(\pa_\cg(v_i)) = \begin{cases} 1 & \text{ if } \pa_\cg(v_i)\text{ contains a vertex of }V\text{ or two vertices of }U, \text{ and}
    \\ 0 & \text{otherwise.}\end{cases}
    \end{equation}
    For $w=u_j$, no edge variable of $\ell$ involves $u_j$, and a pair of parents contributes $1$ if it contains a vertex of $U$ (through the variables $x_{v_iu_ju_k}$ and $x_{u_iu_ju_k}$) and $0$ otherwise. Hence 
    \[
    \ell_{u_j}(\pa_\cg(u_j)) = \begin{cases}
        1 & \text{ if $u_j$ has two parents and at least one of them is in $U$, and} \\
        0 & \text{ otherwise.}
    \end{cases}
    \]

    So every vertex of $W$ contributes $0$ or $1$, and it remains to show that at least two vertices of $W$ do not contribute to $\ell$. A vertex of $W$ without parents in $W$ never contributes. Let $w_1$ be the first vertex of $W$ in a topological order of $\cg$; it does not contribute. If $w_1\in V$, then the first vertex of $U$ in this order has no parent in $U$, so it does not contribute either. If $w_1\in U$, let $w_2$ be the next vertex of $W$. Its only possible parent in $W$ is $w_1\in U$; so if $w_2\in U$ it does not have two parents in $W$, and if $w_2\in V$ it has no parent in $V$ and at most one in $U$. In both cases $w_2$ does not contribute. Hence $\ell(\bc_\cg)\leq p+q-2$.
 
    \textbf{Facet-defining:} Let $F$ be the face of $\bdcim{p+q}{2}$ defined by the inequality. By the validity argument, a DAG on $W$ attains the bound $p+q-2$ exactly when exactly two vertices of $W$ do not contribute to $\ell$. For a DAG $\cg$, a vertex $x$ and a set $P$ of at most two vertices, we write $\cg[x\leftarrow P]$ for the digraph obtained from $\cg$ by replacing the parent set of $x$ by $P$ and leaving all other parent sets unchanged. We note that this operation does not always result in a DAG; however, in the specific cases we apply it to below, we will check that the resulting digraph is indeed acyclic. We consider the matrix $M$ whose columns are the imsets of the following DAGs; see Figure~\ref{fig:clique-isolated-cols} for an example of each kind.
    \begin{enumerate}[label=\textbf{(B\arabic*)}]
        \item For $s\in[q]$ and $t\in[p]$, let $\cb_{s,t}$ be the DAG in which $v_s$ and $u_t$ are sources, every $v_i$ with $i\neq s$ has the unique parent $v_s$, and every $u_j$ with $j\neq t$ has the two parents $v_s,u_t$. Its non-contributing vertices are $v_s$ and $u_t$.
        \item For each $s \in [q]$, $t \in [p]$ and $i\in[q]\setminus\{s\}$, the DAGs $\cb_{s,t}[v_i\leftarrow P]$ with $P$ equal to $\{v_a\}$, $\{v_a,v_b\}$, $\{v_a,u_k\}$ or $\{u_k,u_l\}$, for all distinct $a,b\in[q]\setminus\{i\}$ and all distinct $k,l\in[p]$.
        \item For each $s \in [q]$, $t \in [p]$ and $j\in[p]\setminus\{t\}$, the DAGs $\cb_{s,t}[u_j\leftarrow P]$ with $P$ equal to $\{v_a,u_k\}$ or $\{u_k,u_l\}$, for all $a\in[q]$ and all distinct $k,l\in[p]\setminus\{j\}$.
        \item For each $s \in [q]$ and $ t \in [p]$, the DAGs $\cb_{s,t}[v_s\leftarrow\{u_t\}]$ and $\cb_{s,t}[u_t\leftarrow\{v_a,v_b\}]$ for all $a\ne b$.
        \item For each $s \in [q], t \in [p]$, and $k\neq t$, let $\mathcal{R}_{s,t,k}$ be the DAG in which $u_k$ is a source, $u_t$ has the unique parent $u_k$, $v_s$ and every $u_j$ with $j\neq k,t$ have the two parents $u_k,u_t$, and every $v_i$ with $i\neq s$ has the unique parent $v_s$. We use $\mathcal{R}_{s,t,k}$ and $\mathcal{R}_{s,t,k}[u_t\leftarrow\emptyset]$. Their non-contributing vertices are $u_k$ and $u_t$.
    \end{enumerate}
    In $\cb_{s,t}$ every vertex other than $v_s,u_t$ is a sink, so the modifications in (B2) and (B3) are acyclic, and we can directly check that the DAGs in (B4) and (B5) are acyclic as well. In (B2) and (B3) the vertex whose parent set is modified contributed $1$ to $\ell$ in $\cb_{s,t}$ and still contributes $1$ in the modified graph, by \eqref{eq: clique and independent set contributions} and the description of the contributions of the $u_j$. Similarly, in (B4) the vertices $v_s$ and $u_t$ contributed $0$ to $\ell$ in $\cb_{s,t}$ and still contribute $0$ after these modifications. So, each of these DAGs has exactly two non-contributing vertices and all of  the imsets of these DAGs lie on $F$.

\begin{figure}[ht]
\centering
\setlength{\tabcolsep}{8pt}
\begin{tabular}{ccc}
\begin{tikzpicture}[>=stealth,x=0.95cm,y=1.1cm,every node/.style={circle,draw=cimnode,fill=cimnodefill,inner sep=0.5pt,minimum size=5.2mm,font=\scriptsize}]
\node[coll] (v1) at (0,1) {$v_1$};
\node (v2) at (1,1) {$v_2$};
\node (v3) at (2,1) {$v_3$};
\node[coll] (u1) at (0,0) {$u_1$};
\node (u2) at (1,0) {$u_2$};
\node (u3) at (2,0) {$u_3$};
\draw[->,thick,cimedge] (v1) -- (v2);
\draw[->,thick,cimedge] (v1) to[bend left=35] (v3);
\draw[->,thick,cimedge] (v1) -- (u2);
\draw[->,thick,cimedge] (v1) -- (u3);
\draw[->,thick,cimedge] (u1) -- (u2);
\draw[->,thick,cimedge] (u1) to[bend right=35] (u3);
\end{tikzpicture} &
\begin{tikzpicture}[>=stealth,x=0.95cm,y=1.1cm,every node/.style={circle,draw=cimnode,fill=cimnodefill,inner sep=0.5pt,minimum size=5.2mm,font=\scriptsize}]
\node[coll] (v1) at (0,1) {$v_1$};
\node (v2) at (1,1) {$v_2$};
\node (v3) at (2,1) {$v_3$};
\node[coll] (u1) at (0,0) {$u_1$};
\node (u2) at (1,0) {$u_2$};
\node (u3) at (2,0) {$u_3$};
\draw[->,thick,cimedge] (v1) to[bend left=35] (v3);
\draw[->,thick,cimedge] (v1) -- (u2);
\draw[->,thick,cimedge] (v1) -- (u3);
\draw[->,thick,cimedge] (u1) -- (u2);
\draw[->,thick,cimedge] (u1) to[bend right=35] (u3);
\draw[->,thick,cimhl] (u2) -- (v2);
\draw[->,thick,cimhl] (u3) -- (v2);
\end{tikzpicture} &
\begin{tikzpicture}[>=stealth,x=0.95cm,y=1.1cm,every node/.style={circle,draw=cimnode,fill=cimnodefill,inner sep=0.5pt,minimum size=5.2mm,font=\scriptsize}]
\node (v1) at (0,1) {$v_1$};
\node (v2) at (1,1) {$v_2$};
\node (v3) at (2,1) {$v_3$};
\node[coll] (u1) at (0,0) {$u_1$};
\node[coll] (u2) at (1,0) {$u_2$};
\node (u3) at (2,0) {$u_3$};
\draw[->,thick,cimedge] (v1) -- (v2);
\draw[->,thick,cimedge] (v1) to[bend left=35] (v3);
\draw[->,thick,cimedge] (u2) -- (u1);
\draw[->,thick,cimedge] (u1) -- (v1);
\draw[->,thick,cimedge] (u2) -- (v1);
\draw[->,thick,cimedge] (u1) to[bend right=35] (u3);
\draw[->,thick,cimedge] (u2) -- (u3);
\end{tikzpicture} \\
{\scriptsize $\cb_{1,1}$} & {\scriptsize $\cb_{1,1}[v_2\leftarrow\{u_2,u_3\}]$} & {\scriptsize $\mathcal{R}_{1,1,2}$}
\end{tabular}
\caption{Columns for $q=p=3$: an example of (B1), one of (B2) (the new parents of $v_2$ in purple), and a graph of type (B5). Purple vertices are the two vertices that do not contribute to $\ell$. Deleting the arc $u_2\to u_1$ from $\mathcal{R}_{1,1,2}$ gives $\mathcal{R}_{1,1,2}[u_1\leftarrow\emptyset]$.}
\label{fig:clique-isolated-cols}
\end{figure}
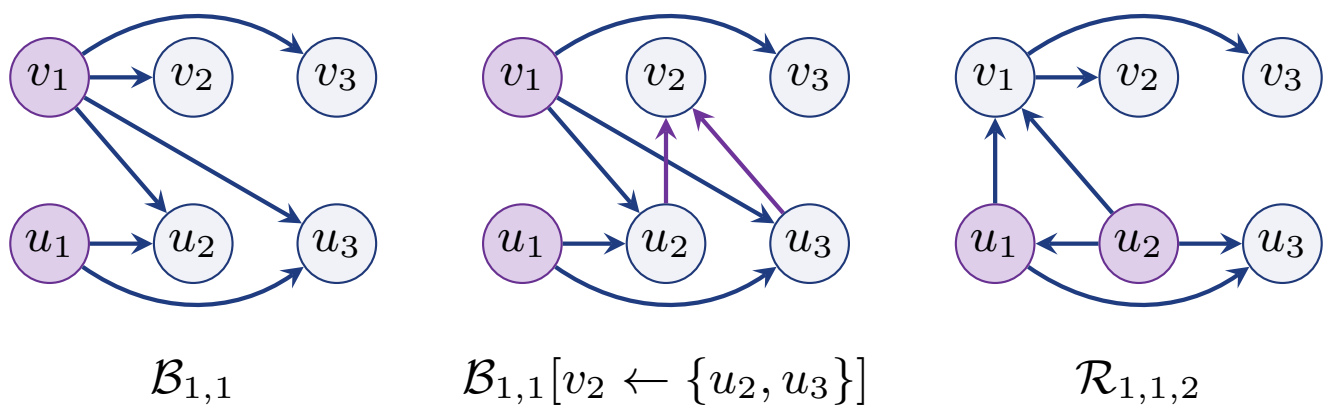

    To show that these columns span $\mathbb R^{\binom{p+q}2+\binom{p+q}3}$, let $f=\sum a_{xy}x_{xy}+\sum b_{xyz}x_{xyz}$ be a linear functional that vanishes on all columns. We will show that $f=0$. By Lemma~\ref{lemma:weighted-coord-count}, $f(\bc_\cg)=\sum_w f_w(\pa_\cg(w))$. So if two columns come from DAGs $\cg$ and $\cg[x\leftarrow P]$, then subtracting them gives $f_x(\pa_\cg(x))=f_x(P)$ since all other parent sets remain unchanged. We apply this to the following pairs.
        \begin{itemize}
        \item \emph{$\cb_{s,t}$ against (B2).} Here $\pa_{\cb_{s,t}}(v_i)=\{v_s\}$, so $f_{v_i}(\{v_s\})=a_{v_sv_i}$. Comparing with $\cb_{s,t}[v_i\leftarrow\{v_a\}]$ gives $a_{v_av_i}=a_{v_sv_i}$ for all $a,s\neq i$. So for each $i$ all $a_{v_av_i}$ are equal, and since $a_{v_av_i}=a_{v_iv_a}$, the coefficients $a_{v_av_b}$ for distinct $a,b\in[q]$ all equal one constant $\kappa$. Comparing $\cb_{s,t}$ with the other DAGs in (B2), whose parent sets of $v_i$ are $\{v_a,v_b\}$, $\{v_a,u_k\}$ and $\{u_k,u_l\}$, then gives
        \[
        b_{v_av_bv_i}=-\kappa,\qquad a_{u_kv_i}+b_{v_av_iu_k}=0,\qquad a_{u_kv_i}+a_{u_lv_i}+b_{v_iu_ku_l}=\kappa
        \]
        for all distinct $a,b,i\in[q]$ and distinct $k,l\in[p]$.
        \item \emph{$\cb_{s,t}$ against (B4).} In $\cb_{s,t}$ the vertex $v_s$ is a source, so comparing with the DAG $\cb_{s,t}[v_s\leftarrow\{u_t\}]$ gives $a_{u_tv_s}=f_{v_s}(\{u_t\})=f_{v_s}(\varnothing)=0$. As $s$ and $t$ are arbitrary, $a_{u_kv_a}=0$ for all $a\in[q]$ and $k\in[p]$, and the previous item then gives $b_{v_av_bu_k}=0$ and $b_{v_au_ku_l}=\kappa$ for all distinct $a,b\in[q]$ and distinct $k,l\in[p]$. The remaining DAGs in (B4) give no new information: in $\cb_{s,t}$ the vertex $u_t$ is a source, so comparing with $\cb_{s,t}[u_t\leftarrow\{v_a,v_b\}]$ gives $a_{v_au_t}+a_{v_bu_t}+b_{v_av_bu_t}=0$, that is, $b_{v_av_bu_t}=0$, which we have already shown.
        \item \emph{$\cb_{s,t}$ against (B3).} Here $\pa_{\cb_{s,t}}(u_j)=\{v_s,u_t\}$, so $f_{u_j}(\{v_s,u_t\})=a_{v_su_j}+a_{u_tu_j}+b_{v_su_tu_j}=a_{u_tu_j}+\kappa$, using $a_{v_su_j}=0$ and $b_{v_su_tu_j}=\kappa$. Comparing $\cb_{s,t}$ with $\cb_{s,t}[u_j\leftarrow\{v_a,u_k\}]$, for which $f_{u_j}(\{v_a,u_k\})=a_{u_ku_j}+\kappa$ in the same way, gives $a_{u_ku_j}=a_{u_tu_j}$. Letting $s$, $t$ and $j$ vary, we see that all $a_{u_ku_l}$ with distinct $k,l\in[p]$ are equal to a constant $\nu$; in particular $f_{u_j}(\{v_s,u_t\})=\nu+\kappa$. Finally, comparing $\cb_{s,t}$ with $\cb_{s,t}[u_j\leftarrow\{u_k,u_l\}]$, for which $f_{u_j}(\{u_k,u_l\})=a_{u_ku_j}+a_{u_lu_j}+b_{u_ju_ku_l}=2\nu+b_{u_ju_ku_l}$, gives $2\nu+b_{u_ju_ku_l}=\nu+\kappa$. Hence $b_{u_ju_ku_l}=\kappa-\nu$ for all distinct $j,k,l\in[p]$.
        \item \emph{$\mathcal{R}_{s,t,k}$ against $\mathcal{R}_{s,t,k}[u_t\leftarrow\emptyset]$.} These differ only in the parent set of $u_t$, which is $\{u_k\}$ or $\emptyset$, so $a_{u_ku_t}=f_{u_t}(\{u_k\})=f_{u_t}(\emptyset)=0$, that is, $\nu=0$.
    \end{itemize}
    Thus, we conclude that $f=\kappa\ell$. But $\ell$ does not vanish on the vertices of $F$, since for instance $\ell(\bc_{\cb_{s,t}})=p+q-2>0$. So we see that $\kappa=0$ and $f=0$. So the columns of $M$ have rank $\binom{p+q}2+\binom{p+q}3$, and thus $F$ is a facet.
\end{proof}

We finally consider the following family of facets whose edge variables correspond to the disjoint union of two cliques. These facets are especially interesting because they are the only ones we have found in which edge variables appear with different signs.

\begin{theorem}\label{thm: two cliques}
    Let $m \geq 3$ and $n \geq 2$. Let $V = \{v_1,\dots,v_m\}$ and $U = \{u_1,\dots, u_n\}$ be disjoint sets of vertices. Then the following inequality defines a facet of $\bdcim{p}{2}$ for $p\geq m+n$:
   
    \begin{equation} \label{eq: two cliques}
    -\sum_{i<j\in[m]} x_{v_iv_j} + \sum_{i<j\in[n]} x_{u_iu_j} + 2\sum_{i<j<k\in[m]}x_{v_iv_jv_k} - \sum_{i<j<k\in[n]}x_{u_iu_ju_k} + \sum_{i<j\in[m],k\in[n]}x_{v_iv_ju_k} \leq n.\end{equation}

\end{theorem}

Before we prove Theorem \ref{thm: two cliques}, let us first consider the following short example for the $m = n = 3$ case. 

\begin{example}\label{eg: two cliques}
    We illustrate Theorem~\ref{thm: two cliques} and its proof in the case $m=n=3$; the proof itself treats general $m$ and $n$. Let $V=\{1,2,3\}$ and $U=\{4,5,6\}$, and we index the variables by $1,\dots,6$, writing $v_i=i$ and $u_k=3+k$. Then the inequality (\ref{eq: two cliques}) reads
    \begin{align*}
        -(x_{12} + x_{13} +x_{23}) + x_{45} + x_{46} +x_{56} + 2x_{123} +x_{124} + x_{125} +x_{126} +x_{134} + x_{135} +x_{136}&\\
        +x_{234} + x_{235} +x_{236}- x_{456}& \leq 3.
    \end{align*}
    Exchanging the roles of $V$ and $U$ gives a second inequality of the same form, in which the edge variables among vertices of $V$ are positive and those among vertices of $U$ are negative.

    \emph{Validity.} Let $\ell$ denote the left-hand side. The proof evaluates $\ell(\bc_\cg)$ one vertex at a time, using $\ell(\bc_\cg)=\sum_w\ell_w(\pa_\cg(w))$ from Lemma~\ref{lemma:weighted-coord-count}. Up to symmetry, the possible parent sets of size at most two give:
    \begin{center}
    \begin{tabular}{c|cccccc}
        $\pa_\cg(w)$ & $\emptyset$ & $\{v\}$ & $\{u\}$ & $\{v,v'\}$ & $\{v,u\}$ & $\{u,u'\}$\\ \hline
        $w\in V$ & $0$ & $-1$ & $0$ & $-1-1+2=0$ & $-1+0+1=0$ & $0$\\
        $w\in U$ & $0$ & $0$ & $1$ & $0+0+1=1$ & $0+1+0=1$ & $1+1-1=1$
    \end{tabular}
    \end{center}
    where $v,v'\in V$ and $u,u'\in U$. So each vertex of $V$ contributes at most $0$ and each vertex of $U$ at most $1$, and hence $\ell(\bc_\cg)\leq 3$. For example, the DAG~(a) in Figure~\ref{fig:two-cliques-cols} has $\pa(u_1)=\pa(u_2)=\{v_1,v_2\}$ and $\pa(u_3)=\{v_2,v_3\}$, so $\ell(\bc_\cg)=x_{124}+x_{125}+x_{236}=3$. On the other hand, the DAG with the single arc $v_1\to v_2$ has $\ell(\bc_\cg)=-x_{12}=-1$.

Reading off the table, $\cg$ attains the bound if and only if each of $u_1,u_2,u_3$ has a parent in $U$ or two parents in $V$, and no vertex of $V$ has a unique parent which lies in $V$. In the proof we construct seven families (a)-(g) of such DAGs, each obtained from a simple base DAG by changing the parent set of a single vertex (Figure~\ref{fig:two-cliques-cols} shows one DAG from each family), and we show that their imsets span all $\binom{m+n}2+\binom{m+n}3$ coordinates. This shows that the inequality defines a facet.
\end{example}

\begin{proof}[Proof of Theorem \ref{thm: two cliques}]

    \textbf{Validity:} Write $V=\{v_1,\dots,v_m\}$, $U=\{u_1,\dots,u_n\}$ and $W=V\cup U$, and let $\ell$ be the left-hand side of the inequality (\ref{eq: two cliques}). Every variable of $\ell$ is indexed by a subset of $W$. For a DAG $\cg\in\dags{p}{2}$, the restriction of $\bc_\cg$ to the coordinates indexed by subsets of $W$ is the imset of the induced subgraph $\cg[W]$, so $\ell(\bc_\cg)=\ell(\bc_{\cg[W]})$. Hence it suffices to prove validity for DAGs on the vertex set $W$. By Lemma~\ref{lem: building facets}, it also suffices to show that the inequality defines a facet of $\bdcim{m+n}{2}$. Throughout the proof we therefore take $p=m+n$ and consider DAGs $\cg$ on the vertex set $W$. By Lemma \ref{lemma:weighted-coord-count}, $\ell(\bc_\cg)=\sum_{w\in W}\ell_w(\pa_\cg(w))$.

    For $w=v_i$, a parent in $V$ contributes $-1$ through an edge variable and a parent in $U$ contributes $0$. A pair of parents contributes $2$ if both are in $V$, $1$ if one is in $V$ and the other in $U$ (through the variable $x_{v_iv_ju_k}$), and $0$ otherwise. Considering all parent sets of size at most two, we see that $\ell_{v_i}(\pa_\cg(v_i))=-1$ if $\pa_\cg(v_i)=\{v_j\}$ for some $v_j\in V$, and $\ell_{v_i}(\pa_\cg(v_i))=0$ otherwise. In particular, $\ell_{v_i}(\pa_\cg(v_i))\leq0$.

    For $w=u_k$, a parent in $U$ contributes $1$ and a parent in $V$ contributes $0$. A pair of parents contributes $1$ if both are in $V$, $-1$ if both are in $U$, and $0$ otherwise. Considering all cases, we see that $\ell_{u_k}(\pa_\cg(u_k))=1$ if $\pa_\cg(u_k)$ contains a vertex of $U$ or is a pair of vertices of $V$, and $0$ otherwise. In particular, $\ell_{u_k}(\pa_\cg(u_k))\leq1$.

    Therefore $\ell(\bc_\cg)\leq\sum_{k\in[n]}\ell_{u_k}(\pa_\cg(u_k))\leq n$, and the inequality is valid on $\bdcim{p}{2}$ for all $p\geq m+n$.

    \textbf{Facet-defining:} Let $F$ be the face of $\bdcim{m+n}{2}$ defined by the inequality (\ref{eq: two cliques}). By the above, the imset of a DAG $\cg$ on $W$ lies on $F$ if and only if
    \begin{itemize}
        \item every $u_k$ has a parent in $U$ or exactly two parents in $V$, and
        \item no $v_i$ has $\pa_\cg(v_i)=\{v_j\}$ for some $v_j\in V$.
    \end{itemize}
    Since $F$ does not contain the origin, it suffices to find DAGs with these properties whose imsets span $\mathbb R^{\binom{m+n}2+\binom{m+n}3}$.

    For a pair $D=\{a,b\}\subseteq[m]$, let $\cb_D$ be the DAG on $W$ whose arcs are $v_a\to u_k$ and $v_b\to u_k$ for all $k\in[n]$. Thus every $u_k$ has the parent set $\{v_a,v_b\}$, and every vertex of $V$ is a source as each element of $V$ other than $v_a$ and $v_b$ is isolated. Let $\cb=\cb_{\{1,2\}}$. For each $c\in[m]$, fix a pair $D_c\subseteq[m]\setminus\{c\}$; this is possible since $m\geq3$. Let $M$ are the imsets of the following DAGs.
    \begin{itemize}
        \item[(a)] $\cb$, and $\cb[u_k\leftarrow\{v_i,v_j\}]$ for each $k\in[n]$ and each distinct $i,j\in[m]$.
        \item[(b)] $\cb[u_k\leftarrow\{u_l\}]$ for each $1\leq l<k\leq n$.
        \item[(c)] $\cb[u_k\leftarrow\{v_i,u_l\}]$ for each $i\in[m]$ and $1\leq l<k\leq n$.
        \item[(d)] (if $n\geq3$) $\cb[u_k\leftarrow\{u_l,u_{l'}\}]$ for each $1\leq l<l'<k\leq n$.
        \item[(e)] $\cb[v_c\leftarrow\{v_a,v_b\}]$ for all $1\leq a<b<c\leq m$.
        \item[(f)] $\cb_{D_c}[v_c\leftarrow\{u_k\}]$ for each $c\in[m]$ and $k\in[n]$.
        \item[(g)] $\cb_{D_c}[v_c\leftarrow\{v_a,u_k\}]$ for each $c\in[m]$, $a\in[m]\setminus\{c\}$ and $k\in[n]$.
    \end{itemize}
    One DAG of each type is drawn in Figure~\ref{fig:two-cliques-cols}.
    
\begin{figure}[ht]
\centering
\setlength{\tabcolsep}{5pt}
\begin{tabular}{cccc}
\begin{tikzpicture}[>=stealth,x=0.95cm,y=0.95cm,every node/.style={circle,draw=cimnode,fill=cimnodefill,inner sep=0.5pt,minimum size=5.2mm,font=\scriptsize}]
\node (n1) at (0,1.2) {$v_1$};
\node (n2) at (1,1.2) {$v_2$};
\node (n3) at (2,1.2) {$v_3$};
\node[coll] (n4) at (0,0) {$u_1$};
\node[coll] (n5) at (1,0) {$u_2$};
\node[coll] (n6) at (2,0) {$u_3$};
\draw[->,thick,cimedge] (n1) -- (n4);
\draw[->,thick,cimedge] (n2) -- (n4);
\draw[->,thick,cimedge] (n1) -- (n5);
\draw[->,thick,cimedge] (n2) -- (n5);
\draw[->,thick,cimedge] (n2) -- (n6);
\draw[->,thick,cimedge] (n3) -- (n6);
\draw[->,thick,cimedge,white] (n4) to[bend right=40] (n6);
\end{tikzpicture} &
\begin{tikzpicture}[>=stealth,x=0.95cm,y=0.95cm,every node/.style={circle,draw=cimnode,fill=cimnodefill,inner sep=0.5pt,minimum size=5.2mm,font=\scriptsize}]
\node (n1) at (0,1.2) {$v_1$};
\node (n2) at (1,1.2) {$v_2$};
\node (n3) at (2,1.2) {$v_3$};
\node[coll] (n4) at (0,0) {$u_1$};
\node[coll] (n5) at (1,0) {$u_2$};
\node (n6) at (2,0) {$u_3$};
\draw[->,thick,cimedge] (n1) -- (n4);
\draw[->,thick,cimedge] (n2) -- (n4);
\draw[->,thick,cimedge] (n1) -- (n5);
\draw[->,thick,cimedge] (n2) -- (n5);
\draw[->,thick,cimedge] (n4) to[bend right=40] (n6);
\end{tikzpicture} &
\begin{tikzpicture}[>=stealth,x=0.95cm,y=0.95cm,every node/.style={circle,draw=cimnode,fill=cimnodefill,inner sep=0.5pt,minimum size=5.2mm,font=\scriptsize}]
\node (n1) at (0,1.2) {$v_1$};
\node (n2) at (1,1.2) {$v_2$};
\node (n3) at (2,1.2) {$v_3$};
\node[coll] (n4) at (0,0) {$u_1$};
\node[coll] (n5) at (1,0) {$u_2$};
\node[coll] (n6) at (2,0) {$u_3$};
\draw[->,thick,cimedge] (n1) -- (n4);
\draw[->,thick,cimedge] (n2) -- (n4);
\draw[->,thick,cimedge] (n1) -- (n5);
\draw[->,thick,cimedge] (n2) -- (n5);
\draw[->,thick,cimedge] (n3) -- (n6);
\draw[->,thick,cimedge] (n4) to[bend right=40] (n6);
\end{tikzpicture} &
\begin{tikzpicture}[>=stealth,x=0.95cm,y=0.95cm,every node/.style={circle,draw=cimnode,fill=cimnodefill,inner sep=0.5pt,minimum size=5.2mm,font=\scriptsize}]
\node (n1) at (0,1.2) {$v_1$};
\node (n2) at (1,1.2) {$v_2$};
\node (n3) at (2,1.2) {$v_3$};
\node[coll] (n4) at (0,0) {$u_1$};
\node[coll] (n5) at (1,0) {$u_2$};
\node[coll] (n6) at (2,0) {$u_3$};
\draw[->,thick,cimedge] (n1) -- (n4);
\draw[->,thick,cimedge] (n2) -- (n4);
\draw[->,thick,cimedge] (n1) -- (n5);
\draw[->,thick,cimedge] (n2) -- (n5);
\draw[->,thick,cimedge] (n4) to[bend right=40] (n6);
\draw[->,thick,cimedge] (n5) -- (n6);
\end{tikzpicture} \\
{\scriptsize (a) $u_3\leftarrow\{v_2,v_3\}$} & {\scriptsize (b) $u_3\leftarrow u_1$} & {\scriptsize (c) $u_3\leftarrow\{v_3,u_1\}$} & {\scriptsize (d) $u_3\leftarrow\{u_1,u_2\}$} \\[4pt]
\end{tabular}\\
\begin{tabular}{cccc}
\begin{tikzpicture}[>=stealth,x=0.95cm,y=0.95cm,every node/.style={circle,draw=cimnode,fill=cimnodefill,inner sep=0.5pt,minimum size=5.2mm,font=\scriptsize}]
\node (n1) at (0,1.2) {$v_1$};
\node (n2) at (1,1.2) {$v_2$};
\node[coll] (n3) at (2,1.2) {$v_3$};
\node[coll] (n4) at (0,0) {$u_1$};
\node[coll] (n5) at (1,0) {$u_2$};
\node[coll] (n6) at (2,0) {$u_3$};
\draw[->,thick,cimedge] (n1) to[bend left=40] (n3);
\draw[->,thick,cimedge] (n2) -- (n3);
\draw[->,thick,cimedge] (n1) -- (n4);
\draw[->,thick,cimedge] (n2) -- (n4);
\draw[->,thick,cimedge] (n1) -- (n5);
\draw[->,thick,cimedge] (n2) -- (n5);
\draw[->,thick,cimedge] (n1) -- (n6);
\draw[->,thick,cimedge] (n2) -- (n6);
\end{tikzpicture} &
\begin{tikzpicture}[>=stealth,x=0.95cm,y=0.95cm,every node/.style={circle,draw=cimnode,fill=cimnodefill,inner sep=0.5pt,minimum size=5.2mm,font=\scriptsize}]
\node (n1) at (0,1.2) {$v_1$};
\node (n2) at (1,1.2) {$v_2$};
\node (n3) at (2,1.2) {$v_3$};
\node[coll] (n4) at (0,0) {$u_1$};
\node[coll] (n5) at (1,0) {$u_2$};
\node[coll] (n6) at (2,0) {$u_3$};
\draw[->,thick,cimedge] (n2) -- (n4);
\draw[->,thick,cimedge] (n3) -- (n4);
\draw[->,thick,cimedge] (n2) -- (n5);
\draw[->,thick,cimedge] (n3) -- (n5);
\draw[->,thick,cimedge] (n2) -- (n6);
\draw[->,thick,cimedge] (n3) -- (n6);
\draw[->,thick,cimedge] (n5) -- (n1);
\end{tikzpicture} &
\begin{tikzpicture}[>=stealth,x=0.95cm,y=0.95cm,every node/.style={circle,draw=cimnode,fill=cimnodefill,inner sep=0.5pt,minimum size=5.2mm,font=\scriptsize}]
\node[coll] (n1) at (0,1.2) {$v_1$};
\node (n2) at (1,1.2) {$v_2$};
\node (n3) at (2,1.2) {$v_3$};
\node[coll] (n4) at (0,0) {$u_1$};
\node[coll] (n5) at (1,0) {$u_2$};
\node[coll] (n6) at (2,0) {$u_3$};
\draw[->,thick,cimedge] (n2) -- (n4);
\draw[->,thick,cimedge] (n3) -- (n4);
\draw[->,thick,cimedge] (n2) -- (n5);
\draw[->,thick,cimedge] (n3) -- (n5);
\draw[->,thick,cimedge] (n2) -- (n6);
\draw[->,thick,cimedge] (n3) -- (n6);
\draw[->,thick,cimedge] (n2) -- (n1);
\draw[->,thick,cimedge] (n5) -- (n1);
\end{tikzpicture} \\
{\scriptsize (e) $v_3\leftarrow\{v_1,v_2\}$} & {\scriptsize (f) $v_1\leftarrow u_2$} & {\scriptsize (g) $v_1\leftarrow\{v_2,u_2\}$} \\[4pt]
\end{tabular}\\
    \caption{One DAG from each of the families (a)-(g) in the proof of Theorem~\ref{thm: two cliques}, for $m=n=3$. Purple vertices have two parents.}
    \label{fig:two-cliques-cols}
    \end{figure}
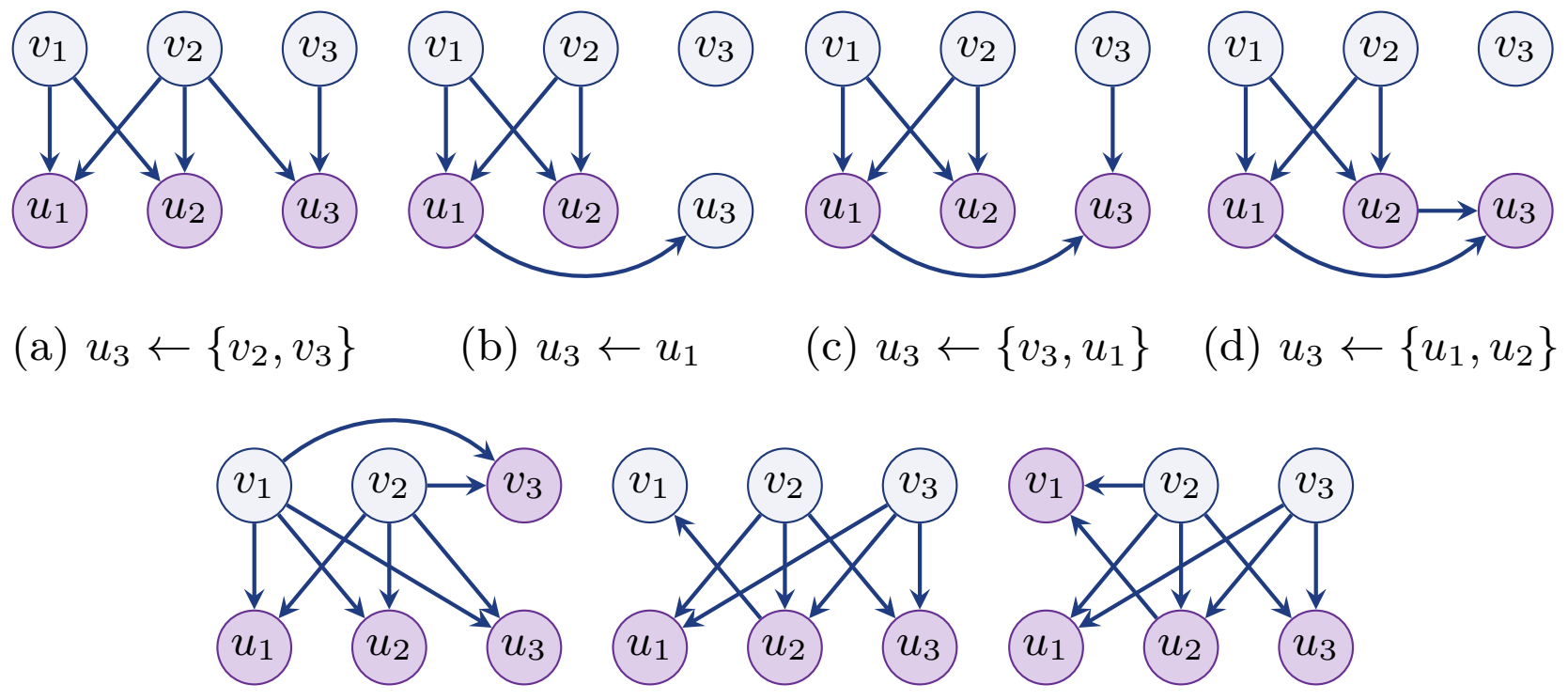

    In $\cb_D$ every vertex of $U$ is a sink. Moreover, every $v_c$ with $c\notin D$ is isolated. So in (a)-(d) we only change the parents of a sink, and in (e)-(g) we only change the parents of an isolated vertex. In both cases no directed cycle is created, so all of these graphs are DAGs. They also lie on $F$. In (a) every $u_k$ has two parents in $V$, and every $v_i$ is a source. In (b)-(d) the modified vertex $u_k$ has a parent in $U$, and all other vertices have the same parent sets as in $\cb$; so every other $u_{k'}$ still has the two parents $v_1,v_2$, and every $v_i$ is still a source. In (e)-(g) every $u_k$ has two parents in $V$, the modified vertex $v_c$ has parent set $\{v_a,v_b\}$, $\{u_k\}$ or $\{v_a,u_k\}$, none of which is a single vertex of $V$, and all other $v_i$ are sources.

    To compute the rank, we use that, by Lemma~\ref{lemma:weighted-coord-count}, modifying the parent set of a single vertex $x$ changes only the coordinates indexed by sets containing $x$:
    \[
    \bc_{\cg[x\leftarrow P]}-\bc_\cg=\be_{x|P}-\be_{x|\pa_\cg(x)},\qquad\text{where } \be_{x|P}=\sum_{y\in P}\be_{xy}+[\,|P|=2\,]\,\be_{P\cup\{x\}} .
    \]
    We split the coordinates into three blocks. The first block consists of the triples of types $v_av_bv_c$, $u_lu_{l'}u_k$ and $v_iu_lu_k$. The second consists of the edges $u_lu_k$. The third consists of the remaining coordinates: the edges $v_av_c$ and $v_iu_k$, and the triples $v_av_bu_k$.
    \begin{itemize}
            \item Each column of (c), (d) and (e) has exactly one nonzero coordinate in the first block, equal to $1$, namely $x_{v_iu_lu_k}$, $x_{u_lu_{l'}u_k}$ and $x_{v_av_bv_c}$ respectively. Every coordinate of the first block arises in this way from exactly one of these columns, so the restrictions of the columns of (c), (d) and (e) to the first block are precisely the unit vectors of that block. All other columns, those of (a), (b), (f) and (g), are zero on the first block.
        \item Among the remaining columns, only those of (b) are nonzero in the second block, and each has exactly one nonzero coordinate there, namely $x_{u_lu_k}$.
    \end{itemize}
    Moreover, the columns of (a), (f) and (g) are zero on the first two blocks: their only modified vertices have parent sets inside $V$ or of the form $\{u_k\}$, $\{v_a,u_k\}$ for a vertex $v_c$, so they only involve the edges $v_av_c$, $v_iu_k$ and the triples $v_av_bu_k$. Order the rows as first, second, third block, and the columns as (c), (d), (e), then (b), then (a), (f), (g). Then
    \[
    M=\begin{pmatrix} P_1&0&0\\ \ast&P_2&0\\ \ast&\ast&N\end{pmatrix},
    \]
    where, after reordering, $P_1$ and $P_2$ are identity matrices of sizes $\binom m3+\binom n3+m\binom n2$ and $\binom n2$, and $N$ is the restriction of the columns of (a), (f) and (g) to the third block. As in the proof of (3), a vanishing linear combination of the columns must have zero coefficients on the columns of $P_1$, by the rows of the first block, then on those of $P_2$, by the rows of the second block, and so we get:
    \[
    \rank(M)\;\geq\;\tbinom m3+\tbinom n3+m\tbinom n2+\tbinom n2+\rank(N).
    \]
    The third block has dimension $\binom m2+mn+n\binom m2$. Let $L$ be the span of the columns from (a), (f) and (g); since these columns are supported on the third block, $\rank(N)=\dim L$, and we show that $L$ has maximal dimension.
    We first collect some vectors in $L$.
    \begin{itemize}
        \item $\bc_\cb\in L$, since $\cb$ is a column of (a). For each $k\in[n]$ and each pair $P\subseteq V$, the column $\cb[u_k\leftarrow P]$ of (a) differs from $\cb$ only in the parent set of $u_k$, so
        \[
        \bc_{\cb[u_k\leftarrow P]}-\bc_\cb=\be_{u_k|P}-\be_{u_k|\{v_1,v_2\}}\in L .
        \]
        \item For a pair $D\subseteq[m]$, the DAG $\cb_D$ is obtained from $\cb$ by changing the parent set of every $u_k$ from $\{v_1,v_2\}$ to $D$, one vertex at a time. Hence $\bc_{\cb_D}-\bc_\cb=\sum_{k\in[n]}\bigl(\be_{u_k|D}-\be_{u_k|\{v_1,v_2\}}\bigr)\in L$, and so $\bc_{\cb_D}\in L$.
        \item In $\cb_{D_c}$ the vertex $v_c$ is isolated, because $c\notin D_c$. So the column $\cb_{D_c}[v_c\leftarrow\{u_k\}]$ of (f) gives $\bc_{\cb_{D_c}[v_c\leftarrow\{u_k\}]}-\bc_{\cb_{D_c}}=\be_{v_c|\{u_k\}}=\be_{v_cu_k}\in L$. Thus $\be_{v_iu_k}\in L$ for all $i\in[m]$ and $k\in[n]$.
        \item The column $\cb_{D_c}[v_c\leftarrow\{v_a,u_k\}]$ of (g) and the column $\cb_{D_c}[v_c\leftarrow\{u_k\}]$ of (f) differ only in the parent set of $v_c$, so
        \[
        g_{ac,k}:=\be_{v_c|\{v_a,u_k\}}-\be_{v_c|\{u_k\}}=\be_{v_av_c}+\be_{v_av_cu_k}\in L
        \]
        for all distinct $a,c\in[m]$ and all $k\in[n]$.
    \end{itemize}
    We now show that every unit vector supported on a row in the third block lies in $L$. The edges $\be_{v_iu_k}$ were obtained above. We now consider the triples.
    \begin{itemize}
        \item By the last item above with $\{a,c\}=\{1,2\}$, $g_{12,k}-g_{12,k'}=\be_{v_1v_2u_k}-\be_{v_1v_2u_{k'}}\in L$ for all $k,k'\in[n]$.
        \item We have $\bc_\cb=\sum_{k\in[n]}\bigl(\be_{v_1u_k}+\be_{v_2u_k}+\be_{v_1v_2u_k}\bigr)$, and the edges $\be_{v_1u_k},\be_{v_2u_k}$ lie in $L$, so $\sum_{k\in[n]}\be_{v_1v_2u_k}\in L$.
        \item Hence $n\,\be_{v_1v_2u_1}=\sum_{k}\be_{v_1v_2u_k}+\sum_k(\be_{v_1v_2u_1}-\be_{v_1v_2u_k})\in L$, so $\be_{v_1v_2u_1}\in L$ as $n\geq 1$, and then $\be_{v_1v_2u_k}=\be_{v_1v_2u_1}-(\be_{v_1v_2u_1}-\be_{v_1v_2u_k})\in L$ for every $k\in[n]$.
        \item For a pair $P=\{v_a,v_b\}$ and $k\in[n]$, the vector from (a) is
        \[
        \be_{u_k|P}-\be_{u_k|\{v_1,v_2\}}=\be_{v_au_k}+\be_{v_bu_k}+\be_{v_av_bu_k}-\be_{v_1u_k}-\be_{v_2u_k}-\be_{v_1v_2u_k} .
        \]
        All terms except $\be_{v_av_bu_k}$ are already known to lie in $L$, so $\be_{v_av_bu_k}\in L$.
    \end{itemize}
    Finally, for distinct $a,c\in[m]$, $\be_{v_av_c}=g_{ac,1}-\be_{v_av_cu_1}\in L$. So $L$ contains every unit vector of the third block, that is, $L$ is the whole third block, and $\rank(N)=\binom m2+mn+n\binom m2$. Since
    \[
    \tbinom m3+\tbinom n3+m\tbinom n2+n\tbinom m2=\tbinom{m+n}3\quad\text{and}\quad\tbinom n2+\tbinom m2+mn=\tbinom{m+n}2,
    \]
    $M$ has full rank, and $F$ is a facet of $\bdcim{m+n}{2}$, and hence of $\bdcim{p}{2}$ for all $p\geq m+n$.
\end{proof}

\section{Outlook}
The degree-bounded characteristic imset polytope $\bdcim{n}{k}$ is a natural subpolytope of $\cim_n$ to consider in the context of causal discovery. In applications, one typically seeks to fit sparse graphs, and the assumption that the in-degree is bounded is commonly used in this area \cite{pearl2009causality}. The dimension of $\bdcim{n}{k}$ grows only polynomially in the number $n$ of nodes of the DAGs involved, which makes linear programs over it much easier to solve, especially for small $k$. In particular, the complexity of one interior point step over $\bdcim{n}{k}$ is
$O\left(m\binom{n}{k+1}^2\right)$, where $m$ is the number of facet-defining inequalities of $\bdcim{n}{k}$. Despite this dimension reduction, our results show that, even though $\bdcim{n}{k}$ has drastically fewer vertices than $\cim_n$, its number of facets still grows exponentially in $n$. Theorem \ref{thm:conn graph} alone gives exponentially many facets, so a straightforward application of interior point methods to solve linear programs over $\bdcim{n}{2}$ will still be exponential in $n$. 

The facets presented here give a complete facet description for $\bdcim{4}{2}$ together with the obvious facets $x_{ijk} \geq 0$ and $\sum_{\{i,j,k\} \subset S} x_{ijk} \leq |S|-2$ for $S \subset [n]$ such that $|S| \geq 3$. However, for $\bdcim{5}{2}$ they are not complete, and computing a full facet description for $\bdcim{5}{2}$ remains out of reach. Nevertheless, we discovered new facets by translating $\bdcim{5}{2}$ such that the origin is in the interior, and then running linear programs over the dual polytope with random objectives. Generically, an optimal solution corresponds to a facet of $\bdcim{5}{2}$. One such additional facet we discovered is
\[
x_{25}+x_{34}+x_{123}+x_{124}+x_{135}+x_{145}\leq 3.
\]
Based on this, we expect that the complexity of not only the number but also the types of facets of $\bdcim{n}{2}$ will continue to grow with $n$, even beyond the immense complexity displayed by the families of facets we have studied here. 

Another advantage of $\bdcim{n}{k}$ that we do not explore here is that it has vastly fewer vertices, and thus the simplex algorithm, whose complexity is determined by the number of vertices and dimension of the ambient space, might be much more tractable for solving linear programs over $\bdcim{n}{k}$ compared to $\cim_n$. In future work it would be interesting to analyze the edges of $\bdcim{n}{k}$ for small $k$ and compare their structure to those which are known for $\cim_n$ \cite{linusson2026edges}. 

\bibliographystyle{plain}
\bibliography{refs}

\end{document}